\documentclass[11pt]{amsart}
\newcommand{\remove}[1]{}

\numberwithin{equation}{section}

\makeindex 

\usepackage{amssymb}

\usepackage{mathrsfs}

\usepackage{graphicx}
\usepackage{amsmath}
\usepackage{amsmath,amssymb}
\usepackage{yfonts}
\usepackage{mathtools}
\usepackage{mathabx}
\usepackage{enumitem}
\usepackage{comment}

 \DeclareMathAlphabet{\mathcalligra}{T1}{calligra}{m}{n}
 
  \DeclareMathOperator{\Aut}{Aut}
 \DeclareMathOperator{\Char}{char}

  \DeclareMathOperator{\depth}{depth}
 \DeclareMathOperator{\End}{End}
 \DeclareMathOperator{\gr}{gr}

 \DeclareMathOperator{\height}{height}

 \DeclareMathOperator{\orbit}{orbit}

 \DeclareMathOperator{\Span}{span}
 \DeclareMathOperator{\spec}{spec}
 \DeclareMathOperator{\rank}{rank}
 
 \DeclareMathOperator{\diag}{diag}
 \DeclareMathOperator{\im}{Im}

\newtheorem{Theorem}{Theorem}[section]
\newtheorem{Theorem--Definition}[Theorem]{Theorem--Definition}
\newtheorem{Corollary}[Theorem]{Corollary}
\newtheorem{Lemma}[Theorem]{Lemma}
\newtheorem{Lemma--Definition}[Theorem]{Lemma--Definition}

\newtheorem *{theorem}{Theorem}
\newtheorem *{theorema}{Theorem A}
\newtheorem *{theoremb}{Theorem B}
\newtheorem *{theoremba}{Theorem B*}
\newtheorem *{theoremc}{Theorem C}
\newtheorem *{theoremd}{Theorem D}
\newtheorem *{theoremda}{Theorem D*}
\newtheorem *{theoreme}{Theorem E}
\newtheorem *{theoremf}{Theorem F}
\newtheorem *{theoremg}{Theorem G}
\newtheorem *{theoremh}{Theorem H}
\newtheorem *{theoremi}{Theorem I}
\newtheorem *{theoremj}{Theorem J}
\newtheorem *{theoremk}{Theorem K}
\newtheorem *{theoreml}{Theorem L}
\newtheorem *{theoremm}{Theorem M}
\newtheorem *{theorem1}{Theorem 1}
\newtheorem *{theorem2}{Theorem 2}
\newtheorem *{question}{Question}
\newtheorem *{definition}{Definition}
\newtheorem *{remark}{Remark}
\newtheorem *{notation}{Notation}
\newtheorem *{note}{Note}
\newtheorem *{conjecture}{Conjecture}
\newtheorem *{conjecturec}{Conjecture C}
\newtheorem *{convention}{Convention}
\newtheorem *{warning}{Warning}
\newtheorem *{proposition}{Proposition}

\newtheorem{Definition-Remark}[Theorem]{Definition-Remark}
\newtheorem{Proposition}[Theorem]{Proposition}
\newtheorem{Proposition--Definition}[Theorem]{Proposition--Definition}

\newtheorem{Note}[Theorem]{Note}

\theoremstyle{remark}

\theoremstyle{definition}
\newtheorem{Remark}[Theorem]{Remark}
\newtheorem{Example}[Theorem]{Example}
\newtheorem{Remark--Definition}[Theorem]{Remark--Definition}
\newtheorem{Definition}[Theorem]{Definition}

\begin{document}
\title{Finite Groups, Smooth Invariants and Isolated Quotient Singularities}
\author{Amiram Braun}
\vspace{.4in}

\date{\today\vspace{-0.5cm}}
\maketitle

\begin{center}
Dept. of Mathematics, University of Haifa, Haifa, Israel 31905

abraun@math.haifa.ac.il
\end{center}

\begin{abstract}
Let $G\subset SL(V)$ be a finite group, $p=\Char F$ and $\dim_FV$ is finite.
Let $S(V)$ be the symmetric algebra of $V$, $S(V)^G$ the subring of $G$-invariants, and $V^*$ the dual space of $V$.
We determine when $S(V)^G$ is a polynomial ring.

\begin{theorema}
$S(V)^G$ is a polynomial ring if and only if:
\begin{enumerate}
\item $G$ is generated by transvections and,
\item $S(V)^{G_U}$ is a polynomial ring for each subspace $U \subset V^*$ with $\dim_F U=1$, where $G_U=\{g \in G|g(u)=u$, $\forall u \in U\}$.
\end{enumerate}
\end{theorema}
As a consequence we obtain the following classification result.
$\hat{A}$ stands for the completion of $A$ at its unique graded maximal ideal.
\begin{theoremb}
Suppose $F$ is algebraically closed and $S(V)^G$ is an isolated singularity.
Then $\widehat{S(V)^G}\cong \widehat{S(W)^H}$, where $\dim_F W=\dim_F V$, $(|H|,p)=1$ and $(H,W)$ is a $\bmod p$ reduction (in the sense of Brauer) of a member in the Zassenhaus-Vincent-Wolf list of complex isolated quotient singularities (e.g. \cite{Wo})
\end{theoremb}


\end{abstract}

\section{\bf Introduction}\label{SEC-Intro}
Let $G\subset GL(V)$ be a finite group, $F$ a field and $\dim_F V$ is finite.
Let $S(V)$ be the symmetric algebra of $V$ and $S(V)^G:=\{x\in S(V)|g(x)=x$, $\forall g \in G\}$.
The basic question we address here is:

\begin{question}
When is $S(V)^G$ a polynomial ring:
\end{question}
This was completely settled in the non-modular case as follows.

\begin{theorem} (Shephard-Todd-Chevalley-Serre, \cite[Sec. 7.2]{Be}).
Suppose $(|G|,p)=1$, where $p=\Char F$ or $p=0$. Then $S(V)^G$ is a polynomial ring if and only if $G$ is generated by pseudo-reflections.
\end{theorem}
In the modular case, that is if $p||G|$, the above question is still open. Special cases are described next.

\begin{theorem} (H. Nakajima, \cite{Na}).
Suppose $G$ is a $p$-group and $F=F_p$, the prime field.
Then $S(V)^G$ is a polynomial ring if and only if $V$ has a basis $\{x_1,...,x_n\}$ such that
\begin{enumerate}
  \item $Fx_1+\cdots +Fx_i$ is a $G$-module for $i=1,...,n$;
  \item $|\orbit (x_1)|\cdots |\orbit(x_n)|=|G|$.
\end{enumerate}
\end{theorem}

\begin{theorem} (Kemper-Malle, \cite{KM}, \cite[p. 119]{DK}).
Suppose $V$ is an irreducible $G$-module. Then $S(V)^G$ is a polynomial ring, if and only if:
\begin{enumerate}
  \item $G$ is generated by pseudo-reflections, and;
  \item $S(V)^{G_U}$ is a polynomial ring, for each subspace $U\subset V^*$, the dual space of $V$, where $G_U=\{g \in G|g(u)=u$, $\forall u \in U\}$ (Steinberg condition).
\end{enumerate}
\end{theorem}
For an arbitrary $V$ we have:

\begin{theorem} (Serre \cite{Se}, \cite[V,$\S$6, Ex. 8]{Bour}).
Suppose $S(V)^G$ is a polynomial ring. Then $(1)$ and $(2)$ hold.
\end{theorem}
The necessity of condition $(2)$ was firstly shown in \cite{St} in case $F=\mathbb{C}$.
So Kemper-Malle's theorem implies that Serre's necessary conditions are also sufficient, if $V$ is also irreducible. This suggests that the irreducibility assumption on $V$ can be removed.

\begin{conjecturec}
$S(V)^G$ is a polynomial ring if and only if:
\begin{enumerate}
\item $G$ is generated by pseudo-reflections, and;
\item $S(V)^{G_U}$ is a polynomial ring for each subspace $U \subset U^*$ with $\dim_FU=1$.
\end{enumerate}
\end{conjecturec}
The shift to $\dim_FU=1$ is a simple consequence of the above theorem by Serre.
A proof of Conjecture C, in case $\dim_F V=3$, is given by D. A. Stepanov \cite{Step2} and Shchigolev-Stepanov \cite{Shst}.

An analog result for $p$-groups was recently established in:

\begin{theorem}(A. Braun \cite{Br}).
Suppose $G$ is a finite $p$-group and $\dim _F V \ge 4$.
Then $S(V)^G$ is a polynomial ring if and only if:
\begin{enumerate}
  \item $S(V)^G$ is a Cohen-Macaulay, and;
  \item $S(V)^{G_U}$ is a polynomial ring, for each $U \subset V^*$ with $\dim_F U=2$.
\end{enumerate}
\end{theorem}
Our main objective here is to affirm the previous conjecture if $G \subset SL(V)$.
We prove the following.

\begin{theorema}
Suppose $G \subset SL(V)$ is a finite group. Then $S(V)^G$ is a polynomial ring if and only if:
\begin{enumerate}
\item $G$ is generated by transvections, and;
\item $S(V)^{G_U}$ is a polynomial ring for each subspace $U \subset V^*$ with $\dim_F U=1$.
\end{enumerate}
\end{theorema}
Recall that if $G\subset SL(V)$, the only possible psudo-reflection $g \in G$ is a transvection, that is $\rank (g-1)=1$ and $(g-1)^2(V)=0$.
In particular $g^p=1$, where $p=\Char F$.

Actually condition $(1)$ can be weakened in Theorem A to:
\begin{enumerate}
\item [(1')] $G$ is generated by its $p$-sylov subgroups.
\end{enumerate}

Indeed let $g \in G$, $g^{p^e}=1$, be one of the generators in (1').
Then $(g-1)^{p^e}=0$ on $V^*$, implying that $\ker_{V^*}(g-1)\neq 0$, and hence $\exists U \subset V^*$, $\dim_FU=1$, with $g \in G_U$.
By (2) $G_U$ is generated by pseudo-reflections which are necessarily transvections, since $G \subset SL(V)$.
So each such $g$ can be replaced by these transvection generators.

The proof of Theorem A is achieved by combining the following two separate theorems.

\begin{theorem1}
Suppose conditions $(1)$ and $(2)$ of Theorem A hold.
Let $W \subset V$ be a proper $G$-submodule.
Set $H:=\{g \in G |(g-1)(V)\subseteq W\}$, $m=\dim_F W$.
Then:
\begin{enumerate}
\item $(WS(V))^H=z_1S(V)^H+\cdots +z_mS(V)^H$, where $z_1,...,z_m \in S(V)^G$  and are part of a minimal homogenous generating set of the polynomial ring  $S(V)^H$;
\item $S(V)^G/(z_1,...,z_m)\cong [S(V)^H/(WS(V))^H]^{G/H}$ and is an isolated singularity as well as Cohen-Macaulay.
\end{enumerate}
\end{theorem1}

\begin{theorem2}
Let $G \subset SL(U)$ be a finite group generated by transvections, where $U$ is an irreducible $G$-module and $F$ is algebraically closed.
Let $A \subset S(U)$ be a graded subring having the following properties:
\begin{enumerate}
  \item $A=F[x_1,...,x_d]$ is a polynomial subring, $\dim A = d =\dim_FU$, where $\{x_1,...,x_d\}$ are homogenous generators;
  \item $G$ acts \underline{faithfully} by graded automorphisms on $A$.
\end{enumerate}
Then $\deg x_1=\cdots =\deg x_d$ and $M:=Fx_1+\cdots +Fx_d$ is an irreducible $G$-module.
\end{theorem2}
Consequently the action of $G$ on $A=S(M)$ is obtained from the linear irreducible action of $G$ on $M$.

The proof of Theorem 2 relies on finite group classification results, listing all the possible pairs $(G,U)$.
One then goes through the list and verifies, case by case, the validity of the statements.
This rather lengthly procedure is dealt with in Section 3.

The proof of Theorem 1 is motivated and partially relies on the practice of changing polynomial ring generators, as initiated in \cite{Br}. The relevance of such a result is suggested by \cite[Lemma 2.13]{Na}.

The proof of Theorem A is then established by choosing $U=V/W$ where $W$ is a maximal $G$-submodule of $V$.
Theorem 2 implies that the action of $G/H$ on the polynomial ring $A:=S(V)^H/(WS(V))^H$ is a \underline{linear action}.
Since $A^{G/H}$ is an isolated singularity, it follows that \cite{KM}  can be applied to conclude that $A^{G/H}$ is a polynomial ring.
However this last step also requires that $G/H$ is generated by pseudo-reflections on $M=Fx_1+\cdots +Fx_d$, (actually on $M^*$), where $A=F[x_1,...,x_d]$ $(d=\dim_F U, \deg x_1=\cdots =\deg x_d)$.
To establish this we crucially use the facts that $G \subseteq SL(V)$ as well as $G/H \subseteq SL(U)$ are generated by transvections (on $V$, respectively $U$), and consequently the (abstract) group $G/H$ is generated by its $p$-sylov subgroups.

The lack of this last fact is a serious obstacle in removing the assumption $G\subseteq SL(V)$ in Theorem A, and thus in proving Conjecture C in full generality.

Theorem B states that a modular isolated quotient singularity, if $G \subset SL(V)$, is isomorphic, after completion to a non-modular one. This in turn can be obtained as the completion of a $\bmod p$ reduction of a complex isolated quotient singularity. The list of all complex isolated quotient singularities is due to Zassenhaus-Vincent-Wolf, and can be seen e.g. in \cite{Wo}. They feature in differential geometry as the solutions to the "Clifford-Klein spherical space form problem".

The key ingredients in the proof of Theorem B are that $S(V)^{T(G)}$ is a polynomial ring, where $T(G):=<g|g$ is a transvection on $V>$, and that $(|G/T(G)|,p)=1$.
Both are consequences of Theorem A.
Now this implies by H. Cartan's theorem that:
\begin{flalign*}
  \widehat{S(V)^G}=\widehat{(S(V)^{T(G)})^{G/T(G)}} \cong \widehat{S(W)^{G/T(G)}},
\end{flalign*}
where $W:=Fy_1+\cdots +Fy_n$, $F[[y_1,...,y_n]]=\widehat{S(V)^{T(G)}}$, and $G/T(G)$ acts linearly on $W$.
Also $W$ is void of fixed points for each $h \in G/T(G)$, $h\neq 1$.
This last fact is a key property in the classification of complex isolated quotient singularities.
The shift, from complex classification results to those over algebraically closed $F$  with $\Char F = p>0$, is possible since $(|G/T(G)|,p)=1$, and this is well known to finite group theorists (and others).
It is sketched (e.g.) in \cite[Thm. 3.13]{Step}.
All of this is detailed in Section 4.

Item (1) of the following is obtained as a consequence of Theorem A (without appealing to Theorem B) and a result of S. Kovacs \cite{Ks}.
\begin{proposition}
Let $G\subset SL(V)$ be a finite group, $\dim _FV$ is finite and $F$ is perfect.
Suppose $S(V)^G$ is an isolated singularity.
Then:
\begin{enumerate}
\item $S(V)^G$ is a rational singularity, and
\item $S(V)^G$ has a non-commutative crepant resolution , provided it is Gorenstein.
\end{enumerate}
\end{proposition}
For the definition of rational singularity in prime characteristic we refer to \cite[Def. 1.3]{Ks}.

\section{\bf Theorem A}\label{SEC-ThmA}
Throughout this section we shall use the following notations:

$W\subset V$ is a $G$-submodule, $H:=\{g\in G|(g-1)(V)\subseteq W\}$.
Clearly $H$ is a normal subgroup of $G$.

\begin{Lemma}\label{ThmA-L1}
Let $y=a_nv^n+\cdots +a_1v+a_0 \in (WS(V))^H$, be a homogenous element with $a_i\in WS(V')$, where $W \subseteq V' \subset V$ is a codimension $1$ subspace of $V$, and $v \in V$, with $Fv+V'=V$.
Then:
\begin{flalign*}
[\sum_{s=0}^{i-1}\phi (a_{n-s})\delta (v)^{i-s}\binom{n-s}{i-s}] +\delta(a_{n-i})=0 \text{, where } \phi \in H \text{ and } \delta :=\phi - 1.
\end{flalign*}
\end{Lemma}

\textbf{Proof:}
$0=\delta(y)=\delta(a_nv^n+\cdots +a_1v+a_0)$.
We shall compute the coefficient of $v^{n-i}$ in the last expression.
Recall that $\delta (v^j)=\sum_{r=0}^{j-1}v^r\delta(v)^{j-r}\binom{j}{r}$ (e.g \cite[p.7]{Br}).
$\delta (a_jv^j)=\delta(a_j)v^j+\phi (a_j)\delta(v^j)=\delta (a_j)v^j+\phi (a_j)(\sum_{r=0}^{j-1}v^r\delta (v)^{j-r}\binom{j}{r})$.

So for $j=n,n-1,...,n-i+1$, $v^{n-i}$ will occur with coefficient $\phi(a_j)\delta (v)^{j-(n-i)}\binom{j}{n-i}$ (taking $r=n-i$).

For $j=n-i$, the coefficient of $v^{n-i}$ is $\delta (a_{n-i})$.
Changing indices by taking $j=n-s$, we get that the coefficient of $v^{n-i}$ in the above expression is:\\
$A_{n-i}:=\sum_{s=0}^{i-1} \phi (a_{n-s})\delta(v)^{i-s}\binom{n-s}{i-s}+\delta(a_{n-i})$, (using  $\binom{n-s}{n-i}=\binom{n-s}{i-s})$.

Now $0=\sum_{i=0}^nA_{n-i}v^{n-i}$.
Also $\delta (V)\subseteq W$ implies that $\phi (a_{n-s})\delta (v)^{i-s}$, $\delta (a_{n-i})\in WS(V')$.
Consequently $A_{n-i}\in WS(V')$ and the above algebraic equation of $v$ over $S(V')$ implies $A_{n-i}=0$, $i=0,...,n$.$\Box$

\begin{Lemma}\label{ThmA-L2}
Keeping the notation of Lemma \ref{ThmA-L1}, we have for $k<n$:
\begin{flalign*}
 \sum_{i=0}^{k}\binom{n-i}{k-i}a_{n-i}v^{k-i}\in (WS(V))^H.
\end{flalign*}
\end{Lemma}
\textbf{Proof:}
Let $\phi \in H $ and $ \delta:=\phi - 1$.
Using $\delta (xy)=\delta (x)y+\phi (x)\delta (y)$, we get:
\begin{flalign*}
 & \delta [\sum_{i=0}^{k}\binom{n-i}{k-i}a_{n-i}v^{k-i}] \\
 & = \sum_{i=0}^{k}\binom{n-i}{k-i} {\delta (a_{n-i})v^{k-i}+\sum_{i=0}^{k-1}\phi(a_{n-i})(\sum_{j=1}^{k-i}v^{k-i-j}\delta (v)^j\binom{k-i}{j}}) \\
 & = \sum_{t=0}^{k}v^{k-t}\{\sum_{s=0}^{t-1}\binom{n-s}{k-s}\phi (a_{n-s})\delta(v)^{t-s}\binom{k-s}{t-s}+\binom{n-t}{k-t}\delta(a_{n-t})\}
\end{flalign*}
(using the convention $\sum^{-1}=0$).
But $\binom{n-s}{k-s}\binom{k-s}{t-s}=\binom{n-s}{n-t}\binom{n-t}{k-t}=\binom{n-s}{t-s}\binom{n-t}{k-t}$, implying that the coefficient of $v^{k-t}$ is:
\begin{flalign*}
& \{[\sum_{s=0}^{t-1}\phi(a_{n-s})\delta (v)^{t-s}\binom{n-s}{t-s}]+\delta (a_{n-t})\} \binom{n-t}{k-t}\\
& =A_{n-t}\binom{n-t}{k-t}=0,
\end{flalign*}
where the last equality is by Lemma \ref{ThmA-L1}.$\Box$

\begin{Proposition}\label{ThmA-P3}
Let $y=a_nv^n+\cdots +a_1v+a_0\in (WS(V))^H$, be an homogenous element, where $W\subseteq V' \subset V$, $Fv+V'=V$, and $a_i\in WS(V')$, $i=0,...,n$.
Then $y=b_nv^n+\cdots +b_1v+a_0$, where $b_i\in (WS(V))^H$, $i=1,...,n$, and $\deg b_i=\deg a_i< \deg y$, $i=1,...,n$.
\end{Proposition}

\textbf{Proof:}
By Lemma \ref{ThmA-L2} (taking $k=n-1$), we get
\begin{flalign*}
b_1=\binom{n}{n-1}a_nv^{n-1}+\binom{n-1}{n-2}a_{n-1}v^{n-2}+\cdots +\binom{2}{1}a_2v+a_1 \in (WS(V))^H.
 \end{flalign*}
Hence:
\begin{flalign*}
y=[-\binom{n}{n-1}+1]a_nv^n+[-\binom{n-1}{n-2}+1]a_{n-1}v^{n-1}+\cdots +[-\binom{2}{1}+1]a_2v^2+b_1v+a_0.
\end{flalign*}
By using Lemma \ref{ThmA-L2} (taking $k=n-2$), we replace $[-\binom{2}{1}+1]a_2$ by:
\begin{flalign*}
b_2:=[-\binom{2}{1}+1]\{\binom{n}{n-2}a_nv^{n-2}+\binom{n-1}{n-3}a_{n-1}v^{n-3}+\cdots +\binom{3}{1}a_3v+a_2\}\in (WS(V))^H
\end{flalign*}
and
\begin{flalign*}
y=&[1-\binom{n}{n-1}-(-\binom{2}{1}+1)\binom{n}{n-2}]a_nv^n+\cdots +\\
  & [1-\binom{3}{2}-(-\binom{2}{1}+1)\binom{3}{1}]a_3v^3+b_2v^2+b_1v+a_0.
\end{flalign*}
One continues in this way, constructing $b_3,...,b_n$.$\Box$

\begin{Note}\label{ThmA-N4} \
\begin{enumerate}
\item $b_1,...,b_n$ depend on the elements $a_1,...,a_n$ and hence on the choices of $V',v$.
\item Proposition \ref{ThmA-P3} is a key result, enabling the replacement of each member of a minimal homogenous generating set of $(WS(V))^H$ by one in $S(V)^G$, as claimed in Theorem 1. The element $y$ is repeatedly rewritten in accordance to the pseudo-reflection used $\psi $.
\end{enumerate}
\end{Note}

We shall also need the following simple observation.

\begin{Lemma}\label{ThmA-L5}
Let $\psi \in G - H$ be a pseudo-reflection. Then $W \subseteq \ker(\psi -1)$.
\end{Lemma}

\textbf{Proof:}
$\psi \notin H$ and hence $(\psi -1)(V)\nsubseteq W$.
If $W \nsubset \ker(\psi -1)$ then $(\psi -1)(W)\neq 0$, and then $(\psi -1)(V)=(\psi - 1)(W)$ (since $\rank (\psi -1)=1)$.
But $(\psi - 1)(W)\subseteq W$ (since $W$ is a $G$-submodule), leading to a contradiction.$\Box$

We shall need the following version of \cite[Prop. 2.13]{Br}.

\begin{Proposition}\label{ThmA-P5.5}
Let $A= \oplus _nA_n$ be an $\mathbb{N}$-graded polynomial ring over a field $F=A_0$ and $\mathfrak{G} \subset \Aut _{\gr} (A)$ a finite subgroup of graded automorphisms.
Let $x \in A^{\mathfrak{G}}$ be an homgenous element.
Set $\mathfrak{H}:=\ker (\mathfrak{G} \rightarrow \mathfrak{G} | _{A/xA})$, $\mathfrak{m} :=A_{+}$, the irrelevant maximal ideal of $A$ and $\mathfrak{n} :=A_{+}^\mathfrak{G}$.
Suppose:
\begin{enumerate}
\item $x \in \mathfrak{m} - \mathfrak{m} ^2$;
\item $A_{\mathfrak{n}} ^{\mathfrak{G}}$ satisfies Serre's $S_{m+1}$ and $R_m$ conditions with $m \ge 2$.
\end{enumerate}
Then $A^{\mathfrak{G}}/(x)=(A^{\mathfrak{H}}/x A^ {\mathfrak{H}}) ^{{\mathfrak{G}}/{\mathfrak{H}}} $,
and $[A^{\mathfrak{G}}/(x)]_{\mathfrak{n}/(x)} $ satisfies $S_m$ and $R_{m-1}$.
\end{Proposition}

\begin{note}
In \cite[Prop 2.13]{Br} the following extra assumption is required: "$\mathfrak{G}/\mathfrak{H}$ acts faithfully on $A^{\mathfrak{H}}/xA^{\mathfrak{H}}$".
This is obsolete in view of \cite[Lemma 2.8]{Na}, by taking $\mathfrak{p}:=xA$, which is $\mathfrak{G}$-stable and $\mathfrak{H}$ is by definition the inertia subgoup of $\mathfrak{p}$.
\end{note}

Actually we shall need the stronger version of Poposition \ref{ThmA-P5.5} where $A_{\mathfrak{n}}^{\mathfrak{G}}$, $[A^{\mathfrak{G}}/(x)]_{\mathfrak{n}/(x)}$ are replaced by $A^{\mathfrak{G}}$, $A^{\mathfrak{G}}/(x)$ (respectively). For that we use the following:

\begin{Lemma}\label{ThmA-L5.6}
Let $B=\oplus_{n \ge 0} B_n$ be a commutative Noetherian graded ring with $B_0=F$, a field.
Let $\mathfrak{n}:= \oplus_{{n > 0}} B_n$.
Then $B_{\mathfrak{n}}$ satisfies $R_m$, $S_m$, if and only if $B$ satisfies $R_m$, $S_m$ (respectively).
\end{Lemma}

\textbf{Proof:}
Let $\mathfrak{p}\in \spec R$ with $\height =m$, $\mathfrak{p}^*:=\text {the largest graded ideal in } \mathfrak{p}$.
Then by \cite[Lemma 1.5.6(a)]{BH} $\mathfrak{p}^*$ is a prime ideal and $\height (\mathfrak{p}) = 1 + \height (\mathfrak{p}^*)$, if $\mathfrak{p}$ is not graded \cite[Thm. 1.5.8(b)]{BH}.
Consequently if $\mathfrak{p}$ is not graded then $\height (\mathfrak{p}^*)=m-1$.
But $\mathfrak{p}^* \subset \mathfrak{n}$ implying, if $B_{\mathfrak{n}}$ satisfies $R_m$, that $B_{\mathfrak{p}^*}$ is regular.
This implies by \cite[Ex. 2.24(a)]{BH} that $B_{\mathfrak{p}}$ is regular.
Hence $B$ satisfies $R_m$.
The converse implication is trivial.

Similarly by \cite[Thm. 1.5.9]{BH} $\depth B_{\mathfrak{p}} = \depth B_{\mathfrak{p}^*} + 1$, if $\mathfrak{p}$ is not graded.
So if $B_{\mathfrak{n}}$ satisfies Serre's $S_m$ condition then $\depth B_{\mathfrak{p}^*} \ge \min (m, \height \mathfrak{p}^*)$.
Hence if $\mathfrak{p}$ is not graded $\depth B_{\mathfrak{p}} = 1+ \depth B_{\mathfrak{p}^*} \ge min (m, 1+\height \mathfrak{p}^*) = min(m, \height \mathfrak{p})$.$\Box$

Using similar arguments we have:

\begin{Lemma}\label{ThmA-L5.7}
Let $B=\oplus_{n \ge 0} B_n$ be a graded Noetherian commutative ring with $B_0=F$, a field and $d = \dim B$.
Suppose $B$ satisfies $R_{d-1}$.
Then $B$ is an isolated singularity.
\end{Lemma}

\textbf{Proof:}
Let $\mathfrak{n}:= \oplus _{n\ge 1} B_n$, the unique graded maximal ideal of $B$.
Let $\mathfrak{p}$ be a maximal ideal  of $B$ which is not graded.
Let $\mathfrak{p}^*$ be the largest graded ideal in $\mathfrak{p}$.
As in Lemma \ref{ThmA-L5.6} $\height (\mathfrak{p})=1 +\height (\mathfrak{p}^*)$.
Hence $\height (\mathfrak{p}^*) \le d-1$ and therefore $B_{\mathfrak{p}^*}$ is regular implying by \cite[Ex. 2.24(a)]{BH} that $B_{\mathfrak{p}}$ is regular.$\Box$

\begin{Corollary}\label{ThmA-C5.8}
Suppose $F$ is perfect.
Then the following are equivalent:
\begin{enumerate}
\item $S(V)^G$ is an isolated singularity;
\item $S(V)^{G_U}$ is a polynoimial ring for each $U \subset V^*$, $\dim _F U=1$.
\end{enumerate}
\end{Corollary}

\textbf{Proof:}
By Lemma \ref{ThmA-L5.7} item (1) is equivalent to: $S(V)^G$ satisfies $R_{d-1}$, where $d=\dim _F V=\dim S(V)^G$.
The rest follows from \cite[Lemma 2.4]{Br}.$\Box$

The following is a consequence of Proposition \ref{ThmA-P5.5} and Lemma \ref{ThmA-L5.6}. Actually it is part of the proof of Theorem 1.

\begin{Lemma}\label{ThmA-L6}
Suppose $H:=\{g\in G|(g-1)(V)\subseteq W\}$ and $(z_1,...,z_k,y_{k+1},...,y_m)=(WS(V))^H$, where $\{z_1,...,z_k,y_{k+1},...,y_m\}$ is part of a minimal homogenous generating set of the polynomial ring $S(V)^H$. Assume also that:
\begin{enumerate}
  \item $\{z_1,...,z_k\}\subseteq S(V)^G$;
  \item $S(V)^G$ satisfies Serre's $R_{d-1}$ condition $(d=\dim _FV)$ and is Cohen-Macaulay.
\end{enumerate}
Then $S(V)^G/(z_1,...,z_k)\cong [S(V)^H/z_1S(V)^H+\cdots +z_kS(V)^H]^{G/H}$, and it satisfies Serre's $R_{d-1-k}$ condition as well as being Cohen-Macaulay.
\end{Lemma}

\textbf{Proof:}
This is done by induction.
The argument for $z_1$, that is:
\\$S(V)^G/(z_1) \cong (S(V)^H/z_1S(V)^H)^{G/H}$ starting with $z_0=0$, is the same as for the following inductive step and is therefore omitted.
By assumption $z_1,...,z_k\in \mathfrak{m} - \mathfrak{m}^2$ where $\mathfrak{m}:=S(V)^H_{+}$, the irrelevant maximal ideal of $S(V)^H$.
Suppose we have by induction that
\begin{flalign*}
S(V)^G/(z_1,...,z_{i-1}\cong [S(V)^H/z_1(S(V)^H+\cdots +z_{i-1}S(V)^H]^{G/H}
\end{flalign*}
satisfies $R_{d-1-(i-1)}$-condition and is Cohen-Macaulay.

Set $I:=z_1S(V)^H+\cdots +z_{i-1}S(V)^H$, $J:= z_1S(V)^H + \cdots + z_iS(V)^H$, $A:= S(V)^H/I$, $\bar{\mathfrak{m}}:=\mathfrak{m}/I$, $\mathfrak{G}:=G/H$, $x:=\bar{z_i} = z_i+(z_1,,...,z_{i-1}) \in S(V)^G/(z_1,...,z_{i-1})$.

$x\in \bar{\mathfrak{m}}-\bar{\mathfrak{m}}^2$, since $\bar{z}_i$ is part of the minimal generating set of the polynomial ring $A$.

Since $\mathfrak{G}=G/H$ acts faithfully on $S(V)^H/(WS(V))^H$ by \cite[Lemma 2.8]{Na} and $J \subset (WS(V))^H$, it follows that $\mathfrak{G}$ acts faithfully on $S(V)^H/J=A/xA$.
Hence $\mathfrak{H}:= \ker(\mathfrak{G} \rightarrow \mathfrak{G}| _{ A/xA})=1$.
Therefore by Proposition \ref{ThmA-P5.5} $A^\mathfrak{G}/(x)\cong (A/xA)^\mathfrak{G}$.
Hence:
\begin{flalign}\label{ThmA-F0.5}
(S(V)^H/I)^{G/H} \text{ }/\text{ } \bar{z_i}(S(V)^H/I)^{G/H}\cong (S(V)^H/J)^{G/H}.
\end{flalign}
By the inductive assumption, the left hand side of (\ref{ThmA-F0.5}) is isomorphic to:
$$[ S(V)^G/(z_1,... ,z_{i-1})] / \bar{z_i}[S(V)^G/(z_1,...,z_{i-1})]=S(V)^G/(z_1,...z_i).$$
In conclusion $S(V)^G/(z_1,....z_i)\cong [S(V)^H/z_1S(V)^H+\cdots +z_iS(V)^H]^{G/H}$.
Also by Proposition \ref{ThmA-P5.5} and Lemma \ref{ThmA-L5.6} $[S(V)^G/(z_1,...,z_i)]$ satisfies $R_{d-1-i}$ and $S_{d-i}$.
Since $\dim S(V)^G/(z_1,...,z_i)=d-i$, it follows that it is Cohen-Macaulay.$\Box$

The following result is classical.
The earliest reference I could find is \cite{Ch}, but it only deals with the local analog.
Similarly \cite[Prop. 2.2.4]{BH}, \cite[Cor. 1.5]{Sa} merely handle the local case.
So here is a direct proof.

\begin{Lemma}\label{ThmA-L7}
Let $A=F[x_1,...,x_n]$ be a graded polynomial ring, where $x_1,...,x_n$ are homogenous elements.
Suppose $A/\mathfrak{n}$ is a graded polynomial ring, where $\mathfrak{n}\subseteq \mathfrak{m}:=(x_1,...x_n)$ is a graded ideal.
Then $\mathfrak{n}=y_{s+1}A+\cdots + y_nA$, $\{y_{s+1},..., y_n\}$ is part of a minimal homogenous generating set of $A$, and $s=\dim A / \mathfrak{n}$
\end{Lemma}

\textbf{Proof:}
We clearly may assume that $s < n$,

Set $\bar{\mathfrak{m}}:= \mathfrak{m}/\mathfrak{n}=(A/\mathfrak{n})_{+}, \bar{A}:=A/\mathfrak{n}$.
Then $s=\dim_{\bar{A}/\bar{\mathfrak{m}}}\bar{\mathfrak{m}}/\bar{\mathfrak{m}}^2 = \dim_{A/\mathfrak{m}}\bar{\mathfrak{m}}/\bar{\mathfrak{m}}^2 = \dim _{A/\mathfrak{m}}\mathfrak{m}/(\mathfrak{m} ^2 +\mathfrak{n})$.
Hence $\{\tilde{x}_1,...,\tilde{x}_n\}$, the images of $\{x_1,...,x_n\}$ in $\mathfrak{m}/(\mathfrak{m} ^2 +\mathfrak{n})$ are linearly dependent.

We may assume that $\{\tilde{x}_1,...,\tilde{x}_s\}$ are linearly independent over $A/\mathfrak{m}$.
Also all the linear dependencies are taking place in the homogenous components of $\mathfrak{m}/(\mathfrak{m}^2+\mathfrak{n})$.
Hence for each $i>s$ we have $I_i\subseteq \{1,...,s\}$ such that $\tilde{x}_i  =\sum_{j\in I_i}\alpha_{ij}\tilde{x}_j$, $\alpha _{ij}\in A$, $\deg \tilde{x}_i=\deg \tilde{x}_j $, $j\in I_i$.
So since $A/\mathfrak{m}=F$ we can take $\alpha_{ij} \in F$.
Therefore $\forall i >s$, $m_i:=x_i -\sum_{j \in I_i}\alpha_{ij}x_j \in \mathfrak{m}^2+\mathfrak{n}$.
But since $\{x_1,...,x_n\}$ are linearly independent $\bmod \mathfrak{m}^2$, it follows that $m_i \notin \mathfrak{m}^2$, $\forall i>s$.

Clearly $\{x_1,...,x_s,m_{s+1},...,m_n\}$ is a minimal generating set of $F[x_1,...,x_n]$.
Let $y_i\in \mathfrak{n}$ be a homogenous element with $y_i-m_i\in \mathfrak{m}^2$, $i=s+1,...,n$.
Then $A=F[x_1,...,x_n]=F[x_1,...,x_s,m_{s+1},...,m_n]=F[x_1,...,x_s,y_{s+1},...,y_n]$.
So $\{y_{s+1},...,y_n\}$ is part of a minimal generating set of $A$ implying that $(y_{s+1})\subset (y_{s+1},y_{s+2})\subset \cdots \subset (y_{s+1},...,y_{n})$ is a proper chain of prime ideals in $A$.
Hence $\height (y_{s+1},...,y_n)=n-s$.
But $(y_{s+1},...,y_n) \subseteq \mathfrak{n}$ and $\height (\mathfrak{n}) = n-\dim A/\mathfrak{n} =n-s$, so $(y_{s+1},...,y_n)=\mathfrak{n}$.$\Box$

\begin{Corollary}\label{ThmA-C8}
Keeping the notation of this section, and assuming that $S(V)^H$ is a polynomial ring, where $H:=\{g\in G|(g-1)(V)\subseteq W\}$.
Then $(WS(V))^H=y_1S(V)^H+\cdots +y_mS(V)^H$, where $m=\dim_FW$ and $\{y_1,...,y_m\}$ is a part of a minimal homogenous generating set of $S(V)^H$.
\end{Corollary}

\textbf{Proof:}
Recall that $(WS(V))^H=WS(V)\cap S(V)^H$.
Then by \cite[Lemma 2.11]{Na} (the start of the proof of item (2), using $G(V,W)=H$), $S(V)^H/(WS(V))^H$ is a polynomial ring.
Consequently by Lemma \ref{ThmA-L7}, with $A:=S(V)^H$, $\mathfrak{n}:=(WS(V))^H$, the result holds.
The last equality follows since $\height ((WS(V))^H)=\height(WS(V))=\dim_FW$.$\Box$

We shall now prove Theorem 1 (of the introduction).

\begin{Theorem}\label{Thm-1}
Suppose $\dim_FV=d$, $G\subset SL(V)$ is generated by transvections and $S(V)^{G_U}$ is a polynomial ring $\forall U \subset V^*$, $\dim _F U=1$.
Let $W \subset V$ be a $G$-submodule and $H:=\{g\in G|(g-1)(V)\subseteq W\}$.
Then:
\begin{enumerate}
\item $(WS(V))^H=z_1S(V)^H+\cdots + z_mS(V)^H$, where $z_i\in S(V)^G$, $i=1,...,m$ and $m=\dim_FW$,
\item $\{z_1,...,z_m\}$ is part of a minimal homogenous generating set of $S(V)^H$;
\item $S(V)^G/(z_1,...,z_m) \cong [S(V)^H/(WS(V))^H]^{G/H}$ and it satisfies Serre's condition $R_{d-1-m}$ as well as being Cohen-Macaulay.
\end{enumerate}
\end{Theorem}

\textbf{Proof:}
It follows from \cite[Lemma 2.4]{Br} that $S(V)^G$ satisfies $R_{d-1}$ and also $S(V)^{G_U}$ is a polynomial ring for each $U\subset V^*$, $\dim_FU\ge 1$.
Consequently, since $H=G_{W^\bot}$ \cite[Lemma 2.1]{Br}, it follows that $S(V)^H$ is a polynomial (where $W^{\bot}=\{f \in V^* |f(W)=0\})$.
Also by \cite[Thm. 3.1]{Ke} it follows that $S(V)^{G_U}$ being a Cohen-Macaulay ring, for each $U\subset V^*$, $\dim_FU\ge 1$ implies that $S(V)^G$ is Cohen-Macaulay.

By Corollary \ref{ThmA-C8} $(W(S(V))^H=y_1S(V)^H+\cdots +y_mS(V)^H$, $m=\dim_FW=\height ((WS(V))^H)$, and $\{y_1,...,y_m\}$ are part of a minimal homogenous generating set of $S(V)^H$.

We order $\{y_1,...,y_m\}$, according to their degree's $\deg y_1\le \deg y_2 \le \cdots \le \deg y_m$.

Observe that if $y_1\in S(W)$, then $y_1\in S(W)\cap S(V)^H=S(W)^H$.
Let $\psi $ be a pseudo-reflection with $\psi \notin H$.
Then by Lemma \ref{ThmA-L5} $W\subseteq \ker (\psi -1)$, hence since $G=<H, \psi $ $|$ $\psi $ is a pseudo-reflection, $\psi \notin H >$, we get that $y_1 \in S(W)^G \subseteq S(V)^G$.

We shall next show that $y_1 \in S(W)^H$.
We assume that $y_1 \notin S(W)$.
Let $W \subseteq V'\subset V$, with $\dim_FV'=\dim_FV-1$, and $v \in V -V'$.
$H$ acts trivially on $V/W$ implying that $V' $ is an $H$-module.
Suppose $y_1 \notin WS(V')$.
Since $y_1\in WS(V)$, then $y_1=\sum_{i=0}^s a_iv^i$, $a_i \in WS(V')$, $i=0,...,s$, $a_s\neq 0$, $s>0$.
It follows that $\deg y_1=\deg a_i+i$, $i=0,...,s$, and $\deg y_1=\deg a_s +s>\deg a_s$.

Let $g \in H$ be arbitrary and $\delta = g-1$. Then:
\begin{flalign}\label{ThmA-F1}
0= \delta (y_1)=\sum_{i=0}^s \delta (a_i)v^i +\sum_{i=0}^sa_i\delta(v^i) +\sum_{i=0}^s\delta (a_i)\delta (v^i).
\end{flalign}
Since $\delta (v^i)=iv^{i-1} \delta (v) + \sum_{j=0}^{i-2}\binom{i}{j}v^j\delta (v)^{i-j}$, equation (\ref{ThmA-F1}) takes the form $0=\delta (a_s)v^s +\sum_{i=0}^{s-1}x_iv^i$, $x_i \in S(V')$.
Since $a_s \in WS(V')$, then $\delta (a_s)\in WS(V')$, so if $\delta (a_s) \neq 0 $ we reach an algebraic dependence of $v$ over $S(V')$.
Therefore $\delta (a_s)=0$, $\forall g \in H$.
So $a_s \in (WS(V'))^H \subseteq (WS(V))^H$ and $\deg (a_s) < \deg y_i$, $\forall i$, an absurd.
So $y_1 \in WS(V')$,
Continuing in this way with $V''\subset V'$, $\dim_FV''=\dim_FV'-1$, $W\subseteq V''$, $v\in V'-V''$, we reach $y_1 \in S(W)^H$, as claimed.

We shall now consider the following properties:
\begin{enumerate}
  \item [(1')]$(WS(V))^H=(z_1,...,z_k,y_{k+1},...,y_m)$ with $\deg z_i=\deg y_i$, $i=1,...,k$, and $\{z_1,...,z_k\} \subseteq S(V)^G$;
  \item [(2')] $\{z_1,...,z_k,y_{k+1},...,y_m\}$ is a part of a minimal homogenous generating set of $S(V)^H$, and;
  \item [(3')]$S(V)^G /(z_1,...,z_k)\cong [S(V)^H/z_1S(V)^H+\cdots +z_kS(V)^H]^{G/H}$, and it satisfies Serre's $R_{d-1-k}$-condition, as well as being Cohen-Macaulay.
\end{enumerate}

It follows from Lemma \ref{ThmA-L6}, Corollary \ref{ThmA-C8} and the previous reasoning that the choice $z_1=y_1$ satisfies $(1')$, $(2')$ and $(3')$.

So we assume that $(1')$, $(2')$ and $(3')$ hold for $k$.
Set $y:=y_{k+1}$, so $\deg y \ge \deg z_i$, $i=1,...,k$.
We shall now proceed to show that we can replace $y$ by $z\in S(V)^G$ and have in place all the requirements of $(1')$, $(2')$ and $(3')$.

Let $\psi \in G-H$, where $\psi $ is a transvection.
Set $V':=\ker (\psi - 1)$, $v\in V-V'$ with $Fv+V'=V$.
By Lemma \ref{ThmA-L5} $W \subseteq V'$.
Also by Proposition \ref{ThmA-P3} we have $y=b_nv^n+\cdots +b_1v+a_0$, $b_i\in (WS(V))^H$, $\deg b_i=\deg a_i$. $i=1,...,n$, $a_0 \in WS(V')$ (as well as $a_i\in WS(V')$, $i=1,...,n$).
Also $\deg b_i =\deg a_i < \deg y$, implying by $(1')$ that $b_i \in z_1S(V)^H+\cdots + z_kS(V)^H$.
Since $z_i \in S(V)^G$, $i=1,...,k$, it follows that $\psi (b_i)$, $\delta (b_i)\in z_1S(V)+\cdots +z_kS(V)$ where $\delta := \psi - 1$.
Also $\delta (a_0)=0$, since $\delta (V')=0$.
Consequently:
\begin{flalign*}
\delta (y)=\sum_{i=1}^n \delta(b_iv^i) +\delta (a_0)=\sum_{i=1}^n (\delta(b_i)v^i+\psi (b_i)\delta (v^i))\in z_1S(V)+\cdots +z_kS(V).
\end{flalign*}
Also $S(V)^H$ is $G$-stable since $H$ is normal in $G$.
So $\delta (y)\in S(V)^H\cap (z_1S(V)+\cdots +z_kS(V))= S(V)^H \cap [z_1S(V)^H+\cdots +z_kS(V)^H]S(V)$.
Since $S(V)^H$ is a polynomial ring and $S(V)$ is Cohen-Macaulay (being a polynomial ring), it follows from \cite[Prop. 1.5.15 and 2.2.11]{BH} that $S(V)$ is projective and hence free over $S(V)^H$.
Hence by \cite[Thm. 7.5(ii)]{Ma} $S(V)^H \cap [z_1S(V)^H+\cdots +z_kS(V)^H]S(V)=z_1S(V)^H+\cdots +z_kS(V)^H$.
Consequently $\delta(y)\in z_1S(V)^H+\cdots +z_kS(V)^H$, for each such $\psi $.
Observe that there is no reference here to $V'$ any longer.
Hence: $\bar{y}:=y+z_1S(V)^H+\cdots +z_kS(V)^H$ satisfies:
\begin{flalign*}
 \bar{y}\in [S(V)^H/z_1S(V)^H +\cdots +z_kS(V)^H]^{G/H}=S(V)^G/(z_1,...,z_k).
\end{flalign*}
Therefore there exists $z\in S(V)^G$, $\deg z=\deg y$ with $z-y\in z_1S(V)^H+\cdots +z_kS(V)^H$.
Clearly $z_i$, $y\in (WS(V))^H$, $i=1,...,k$, imply $z\in (WS(V))^H \cap S(V)^G$.
Also $z_1S(V)^H+\cdots +z_kS(V)^H+yS(V)^H=z_1S(V)^H+\cdots +z_kS(V)^H+zS(V)^H$.
Moreover if $z-y\neq 0$, then $\deg (z-y)=\deg y \ge \deg z_i$, $i=1,...,k$, so $z-y=\sum_{i=1}^k z_ic_i$, $c_i \in S(V)^H$ and $c_i$ when expressed in the homogenous generators of $S(V)^H$, the ones which appear must have degree smaller than $\deg y$, so $y$ does not appear in the expression of $c_i$ in these generators.
So $z$ can replace $y$ in the minimal homogenous generating set of $S(V)^H$.
All in all we verified with $z_{k+1}:=z\in S(V)^G$, that:
\begin{enumerate}
  \item [(1')] $(WS(V))^H =z_1S(V)^H+\cdots +z_{k+1}S(V)^H+y_{k+2}S(V)^H+\cdots + y_mS(V)^H$;
  \item [(2')] $\{z_1,...,z_k,z_{k+1},y_{k+2},...,y_m\}$ is part of a minimal homogenous generating set of $S(V)^H$;
\end{enumerate}
Finally (3') follows from Lemma \ref{ThmA-L6}.

The theorem now follows by taking $k=m$.$\Box$

\underline{The proof of Theorem A}

We are given that $G \subset SL(V)$ is generated by transvections on $V$ (or $V^*$) and $S(V)^{G_U}$ is a polynomial ring for each $U \subset V^*$ with $\dim_F U=1$.
We want to show that $S(V)^G$ is a polynomial ring.

As in the proof of \cite[Thm. B, p.20]{Br} we may assume that $F$ is algebraically closed, as well as $S(V)^G$ satisfies $R_{d-1}$, where $d=\dim_F V$ and $S(V)^G$ is Cohen-Macaulay.

Let $W \subset V$ be a maximal $G$-submodule, that is $0 \neq V/W$ is a simple $G$-module. Set $H:=\{g\in G|(g-1)(V)\subseteq W\}$. Hence by \cite[Lemma 2.1]{Br} $H=G_{W^\bot}$, where $W^{\bot}:=\{f\in V^* | f(W)=0\}$, implying that $S(V)^H$ is a polynomial ring.

By Theorem \ref{Thm-1} (= Theorem 1) $(WS(V))^H=z_1S(V)^H+\cdots +z_mS(V)^H$, with $\{z_1,...,z_m\}$ are in $S(V)^G$, $m=\dim_FW=\height (WS(V))^H$.
Moreover $z_1,...,z_m$ are part of a minimal homogenous generating set of $S(V)^H$, and $S(V)^G/(z_1,...,z_m)\cong [(S(V)^H/(WS(V))^H)]^{G/H}$, satisfying $R_{d-1-m}$ and the Cohen-Macaulay property.

As in Corollary \ref{ThmA-C8} $A:=S(V)^H/(WS(V))^H$ is a polynomial ring with $\dim A=d-m=\dim_F V/W$. By \cite[Lemma 2.8]{Na} $G/H$ acts faithfully on $A$, this action is the restriction of the faithful (linear and irreducible) action of $G/H$ on $V/W$ and on $S(V/W)$.
By \cite[Lemma 2.6]{Br} $G/H$ is generated by the transvections $\{\sigma _1,..., \sigma _t \}$ on $V/W$ (and on $(V/W)^*$).

Let $A=F[x_1,...,x_{d-m}]$, where $x_i$ is homogenous, $i=1,...,d-m$.
By Theorem \ref{Thm2-T1} (=Theorem 2, which is proved in section 3) $\deg x_1=\cdots =\deg x_{d-m}$ and $M=Fx_1+\cdots +Fx_{d-m}$ is an irreducible $G/H$-module.
Consequently $A=S(M)$ and the action of $G/H$ on $A$ is induced from a linear faithful action of $G/H$ on $M$.
Since $S(M)^{G/H}=A^{G/H}$ satisfies $R_{d-1-m}$ and $F$ is algebraically closed, it follows from \cite[Lemma 2.4]{Br} that $S(M)^{(G/H)_U}$ is a polynomial ring for each subspace $U \subset M^*$ with $\dim_FU=1$.

We have that $\sigma _i^p=1$, $i=1,...,t$.
So $\sigma _i^p=1$ also holds with $\sigma _i$, regarded now as an automorphism on $A=S(M)$ (although it may fail to act as a pseudo-reflection on $M^*$).
However $(\sigma _i-1)^p=0$ on $M^*$, showing that $\ker(\sigma _i -1)\neq 0$ (on $M^*$), hence $\exists U_i \subset M^*$, $\dim _F U_i=1$ with $\sigma _i \in (G/H)_{U_i}$, $i=1,...,t$.
Since $S(M)^{(G/H)_{U_i}}$ is a polynomial ring it follows that $(G/H)_{U_i}$ is generated by the pseudo-reflections on $M^*$, $\forall i$, implying that $G/H$ is generated by pseudo-reflections on $M^*$.

To conclude, $G/H$ acting faithfully, linearly and irreducibly on $M$ and hence on $M^*$, is generated by pseudo-reflections on $M^*$ and $S(M)^{(G/H)_U}$ is a polynomial ring $\forall U \subset M^*$ with $\dim_FU=1$.
Invoking \cite[Main theorem]{KM} we get that $A^{G/H}\cong S(V)^G/(z_1,...,z_m)$ is a polynomial ring in $d-m$ generators.
Hence $S(V)^G$ is generated by $m+(d-m)=d$ elements, showing that it is a polynomial ring.$\Box$

Theorem 1 and Theorem 2 imply the following necessary conditions.

\begin{Proposition}\label{ThmA-P11}
Let $G \subset SL(V)$, $p=\Char F$ and $S(V)^G$ is a polynomial ring.
Let $W \subset V$ be a maximal $G$-submodule.
Set $H=\{ g \in G|(g-1)(V)\subseteq W \}$.
Then one of the following holds:
\begin{enumerate}
\item $\dim _F V/W =1$;
\item $\dim _F V/W=2$, $G/H \cong SL_2(q)$, $p|q$, or $G/H=SL_2(5)$, $p=3$;
\item $n:=\dim_F V/W \ge 3$, $G/H \cong SL_n(q)$, $p|q$;
\item $n:=\dim_F V/W = 3$, $|G/H|=3.A_6$;
\item $n:=\dim_F V/W = 4$, $G/H=SO_4^-(q)$, $q$ even, $p=2$;
\item $p=2$, $G/H \cong D \rtimes S_n$, where $n:=\dim _F V/W> 1$, $D=\{(\xi ^{a_1},...,\xi ^{a_n})|\sum _{i=1}^n a_i \equiv 0(m) \}$, where $\xi$ is a primitive $m$-th root of unity.
\end{enumerate}
\end{Proposition}

\textbf{Proof:}
By Theorems 1 and 2, $S(M)^{G/H}$ is a polynomial ring, where $A:=S(V)^H/(WS(V))^H=F[x_1,...,x_n]=S(M)$, $M:=Fx_1+ \cdots + Fx_n$, $\deg x_1= \cdots =\deg x_n$ and $M$ is a faithful irreducible $G/H$-module.
So we are in the setup analyzed by \cite[Thm. 7.2]{KM}, so the above possibilities are the only ones admitting a polynomial ring of invariants. $\Box$

\section{\bf Theorem 2}\label{Sec-Thm2}
In this section we shall prove the following (Theorem 2 of the introduction).

\begin{Theorem} \label{Thm2-T1}
  Let $F$ be an algebraically closed field with $\Char F=p>0$.
  Let $G \subset GL(U)$ be a finite group generated by transvections, where $U$ is an irreducible finite dimensional $G$-module.
  Let $A\subset S(U)$ be a graded polynomial subring , having the following properties:

  \begin{enumerate}
  \item $A=F[x_1,...,x_n]$ is a polynomial subring, $\dim A = n =\dim_FU$, where $\{x_1,...,x_n\}$ are homogenous generators;
  \item $G$ acts \underline{faithfully} by graded automorphisms on $A$.
\end{enumerate}
Then $\deg x_1=\cdots =\deg x_n$ and $F x_1+\cdots Fx_n$ is an irreducible $G$-module.
\end{Theorem}

\begin{Remark}\label{Thm2-R2}
  The possibilities for $G$ as defined above are naturally divided into two cases:
  \begin{enumerate}[label=(\roman*)]
  \item $G$ acts primitively on $U$, or
  \item $G$ acts imprimitively on $U$.
  \end{enumerate}
  In case (i) the possibilities for $G$ are listed (following \cite[Thm. 2]{Ka}, \cite{ZS}) in \cite[Thm. 1.5]{KM} and are as follows:

\begin{enumerate}[label=(\alph*)]
    \item $G=SL_n(q)$, $Sp_n(q)$, $SU_n(q)$, with $(n,q)\neq (3,2)$, $p|q$;
    \item $n \ge 4$ is even, $p=2$, $G=SO_n ^{\pm} (q)$ with $q$ even;
    \item $n \ge 6$ is even, $p=2$, $G=S_{n+1}$ or $G=S_{n+2}$;
    \item $(n,p)=(2,3)$, $G=SL_2(5)$;
    \item $(n,p)=(3,2)$, $G=3\cdot A_6$,
    \item $(n,p)=(6,2)$, $G=3\cdot U_4(3)\cdot 2$.
  \end{enumerate}
  In case (ii) \cite [1.8, 1.9]{ZS} imply that $p=2$ and $G$ is a monomial subgroup.

  We shall deal with cases (i) and (ii) separately.

  It is of interest to remark that the assumption $p \ge 3$ simplifies the proof considerably as follows:
  \begin{itemize}
  \item {if case (i) holds, we merely need to handle cases (a) and (d);}
  \item {Case (ii) does not occur.}
  \end{itemize}
\end{Remark}

\underline{The proof of Theorem 2 in case (i)}

We assume that $G$ acts primitively on $U$.

Suppose we have, by negation, $\deg x_1= \cdots = \deg_{x_{r_1}} =a_1 < \deg x_{r_1+1}=\cdots = \deg x_{r_2}=a_2 < \cdots < \deg x_{r_{k-1}+1}= \cdots = \deg x_{r_k}=a_k$, $r_k=n$.

Set $ M_1:= Fx_1+\cdots +Fx_{r_1}$, $M_i:=(B_i+Fx_{r_{i-1}+1}+ \cdots +Fx_{r_i})/B_i$, where $B_i=F[x_1, \cdots , x_{r_{i-1}}]_{a_i}$, the homogenous subspace of degree $a_i$, $i=2,...,k$. Clearly $M_i$ is a $G$-module, with $\dim _F M_i <n$ for each $i$.  In fact $\dim _FM_i =r_i- r_{i-1}$, where $r_0=0$, and $\sum_{i=1}^k \dim _F M_i=n$.

Let $0=X_{0i} \subset X_{1i} \subset \cdots X_{ji} \cdots \subset X_{li}=M_i$, be a $G$-submodule decomposition series of $M_i$ $\forall i$. That is $X_{ji} /X_{j-1,i}$ is an irreducible $G$-module $\forall i,j$.
Therefore $\dim_F X_{ji}/X_{j-1,i} <n$, $\forall i,j$.

In each of the possible cases for $G$, we find, in Corollaries \ref{Thm2-C5}, \ref{Thm2-C7}, \ref{Thm2-C9}, \ref{Thm2-C13}, Lemma \ref{Thm2-L11}, and Propositions \ref{Thm2-P11-2}, \ref{Thm2-P11-5}, subgroups $L \subseteq H \subseteq G$ with the following properties:
\begin{enumerate}[label=(\Roman*)]
\item $L$ is not a $p$-group;
\item $Y_k/Y_{k-1}$ is a trivial $L$-module, $1 \le k \le m$, for any decomposition series of $H$-modules of $X_{ji}/X_{j,i-1}$, $\forall i,j$:
$$0=Y_0 \subset Y_1 \subset \cdots \subset Y_m = X_{ji} /X_{j, i-1}.$$
\end{enumerate}

Consequently $\exists t$ such that $(g-1)^{p^t}(M_i)=0$, $i=1,...,k$, $\forall g \in L$.
Therefore $\exists s$, such that $(g-1)^{p^s}(x_j)=0$, $\forall j=1,...,n$, $ \forall g \in L$.
This shows by assumption (2) (the faithful action of $G$ on $A$) that $L$ is a $p$-group, a contradiction.
Consequently $\deg x_1=\cdots =\deg x_n$.

A similar argument shows that $M_1=Fx_1+\cdots +Fx_n$ is an irreducible $G$-module.

This ends the proof of case (i), in the subcases (a), (b), (c) and (f).

Case (d) is handled in Lemma \ref{Thm2-L14} and case (e) in Corollary \ref{Thm2-C16}.$\Box$

We recall (e.g \cite[Definition (5.3.1)]{KL}) the following.

\begin{Definition} \label{Thm2-D3}
Let $G$ be a finite group and $F$ a field. We define:

$R_F(G)= min\{m | \exists$ injective homomorphism $\varphi  : G \rightarrow PGL_m(F)\}$,

$R_p(G)= min\{R_F(G) | F$ is a field, $\Char F=p\}$,

$R_{p'}(G)=min\{R_s(G)| s \text{ prime, } s\neq p\}$.
\end{Definition}

\begin{Lemma}\label{Thm2-L4}
Let $F$ be an algebraically closed field, with $\Char F=p>0$.
Let $H$ be a finite group with $H/Z(H)$ being simple non-abelian.
Assume that $(|Z(H)|,p)=1$ and $p$ divides $|H|$.
Let $M$ be a non-trivial irreducible $FH$-module. Then one of the following holds:
\begin{enumerate}
\item $\dim_FM \ge R_p(H/Z(H))$;
\item $M$ is a trivial $H_p$-module, where $H_p=<P|P$ is a $p$-sylov subgroup of $H>$.
\end{enumerate}
\end{Lemma}

\textbf{Proof:}
Set $\varphi:H\rightarrow GL(M)$, the natural homomorphism and $N:=\ker \varphi$.

Suppose firstly that $N \not \subset Z(H)$.
Then $NZ(H)/Z(H)$ is a non-trivial normal subgroup of $H/Z(H)$, so $NZ(H)=H$.
Since by assumption $N \neq H$, it follows that $1 \neq H/N=NZ(H)/N \cong Z(H)/N \cap Z(H)$, which is by assumption, a group of order prime to $p$.
Therefore $H_pN/N=N/N$ and $H_p\subseteq N$ as claimed.

Suppose now that $N \subseteq Z(H)$.
Assume firstly that $\varphi (H) \subseteq Z(GL(M))=F1_M$.
$H_p$ is generated by elements of $p$-power order $\{g_1,...,g_r\}$.
Hence $\varphi (g_i)^{p^{e_i}}=1_M$.
But $\varphi (g_i)=\lambda_i 1_M$, $\lambda _i \in F$, $i=1,...,r$.
Hence $\lambda _i^{p^e_i}=1$, and $\lambda _i=1$, $i=1,...,r$.
Hence $\varphi (g_i)=1_M$, $i=1,...,r$ and $H_p \subseteq N$. as claimed.

Finally assume that $\varphi (H) \not \subset Z(GL(M))$.
By Schur's Lemma $\varphi(Z(H))\subseteq F1_M=Z(GL(M))$.
Therefore the induced map $\bar{\varphi} :H/Z(H) \rightarrow GL(M)/Z(GL(M))=PGL(M)$ is a non-trivial injective homomorphism.
Consequently $\dim_FM \ge R_p(H/Z(H))$.$\Box$

\begin{Corollary}\label{Thm2-C5}
Let $F=\bar{F}$ and $\Char F = p$.
Let $H=SL_n(q)$, $Sp_n(q)$ or $SU_n(q)$, where $p|q$ and $(n,p)\neq (3,2)$.
Let $M$ be an irreducible $FH$-module.
Then either  $\dim _F M \ge n$ or $M$ is a trivial $H_p=H$-module.
\end{Corollary}

\textbf{Proof:}
We have (respectively) $|Z(H)|=(n,q-1)$, $(2, q-1)$ or $(n,q+1)$.
Hence $(p, |Z(H)|)=1$.
Now $H/Z(H)$ is simple (unless $H=SL_2(2)$, or $H=SL_2(3))$.
So all requirements of Lemma \ref{Thm2-L4}  are in place.
Moreover $H$ is generated by elements of $p$-order (they are transvections in the natural representation of $H$).
Consequently $H=H_p$.
Also $R_p (H/Z(H))=n$, where the last equality is by \cite[Tab. 5.4.C]{KL}.
Suppose now that $H=SL_2(2)$ or $SL_2(3)$ and $\dim M=1$, then since $H$ is generated by transvections, it follows that $M$ is a trivial $H$-module.$\Box$

\begin{Definition}\label{Thm2-D6}
The complex reflection Mitchell group $G_{34}$ has the structure $G_{34}=6.PSU_4(3).2$ (Wikipedia).
In fact $Z(G_{34})=<h>$, where $h=\diag (\xi,\xi,\xi,\xi,\xi,\xi)$,  $\xi^6=1$, $\xi \in \mathbb{C}$ (written in its $6$-dimensional natural  representation).
$W_2(G_{34})$ denotes the  $\bmod 2$  reduction of $G_{34}$.
That is $G_{34}/{<\diag(-1,-1,-1,-1,-1,-1)>}=W_2(G_{34})$.
If $\bar{h}$ is the image of $h$ in $W_2(G_{34})$, then $\bar{h}=\diag(\bar{\xi},\bar{\xi},\bar{\xi}, \bar{\xi}, \bar{\xi}, \bar{\xi} )$, $\bar{\xi} ^3=1$, $\bar{\xi} \in F (\Char F=2)$.
Moreover $Z(W_2(G_{34}))=<\bar{h}>$ and $W_2(G_{34})$ has the structure $W_2(G_{34})=3.U_4(3).2$.

Let $\pi : W_2(G_{34}) \rightarrow W_2(G_{34})/<\bar{h}>$ be the natural projection.
Set $H:=\pi ^{-1}(U_4(3))$.
We clearly have $H/<\bar{h}>\cong U_4(3)$, a simple group, implying that $Z(H)=<\bar{h}>=Z(W_2(G_{34}))$.
In particular $(|Z(H)|,2)=1$, but $2||H|$, and $H_2$ is not a $2$-group.
\end{Definition}

\begin{Corollary}\label{Thm2-C7}
Let $G=W_2(G_{34})$ and $H \subset G$ as in Definition \ref{Thm2-D6}.
Let $M$ be a non-trivial irreducible $FH$-module, where $F=\bar{F}$ and $\Char F=2$.
Then either $\dim_F M \ge 6$ or $M$ is a trivial $H_2$-module.
\end{Corollary}

\textbf{Proof:}
By Definition \ref{Thm2-D6} all the requirements of Lemma \ref{Thm2-L4} are valid for $H$ and $H/Z(H)\cong U_4(3)$.
Now by \cite[Thm 5.3.9, Tab. 5.3.A]{KL} we have $R_2(U_4(3))\ge R_{3'}(U_4(3))=6$.$\Box$

\begin{Definition}\label{Thm2-D8}
Let $G=SO^{\pm}_n(q)$, $n\ge8$ is even,  $q$ is even.
Set $H:=\Omega ^\pm_n(q)$.
Then $PH $ is a simple group \cite[Tab. 5.1.A]{KL}.
Moreover $Z(H)=1$ by \cite[Prop. 2.9.3]{KL}.
Hence $PH=H$.
\end{Definition}

\begin{Corollary}\label{Thm2-C9}
Let $G$, $H$ be as in Definition \ref{Thm2-D8}.
Let $M$ be a non-trivial irreducible $FH$-module, where $\Char F=2$, $F=\bar{F}$.
Then $\dim_F M \ge n$.
\end{Corollary}

\textbf{Proof:}
By \cite[Prop. 5.4.13]{KL} we have $R_2 (H)=n$.
By the simplicity of $H$ we have that $H\subset GL(M)$ and consequently $H \subset PGL(M)$.
Therefore $\dim _F M \ge n$.$\Box$

The missing cases: $SO^{\pm}_4(q)$ and $SO_6^{\pm} (q)$, where $q$ is even, will be treated next.

\begin{Note}\label{Thm2-N10} \underline{The case of $SO^+_4(q)$, $q$ is even.}

Using \cite[(2.5.11), (2.5.12), Section 2.5, Description 4, Prop. 2.9.1(iv), Prop. 2.9.3]{KL} and the notations there, we have: $I=O^+_4(q)=SO^+_4(q)=S$, $\Omega =\Omega ^+_4(q)$, $S=\Omega<g>$, where $g$ is any reflection.

In fact $\Omega .2=S$, where in case of $q>2$, $\Omega $ is the unique subgroup of index $2$. Actually $\Omega = S'$.
In this case we also have that:
$\Omega ^+_4(q)=L_1L_2$,
$L_i \lhd \Omega^+_4(q)$, $i=1,2$, $L_1\cap L_2=1$, $L_i\cong SL_2(q)$, $i=1,2$ (and the latter is a simple group since $q>2$).
Moreover $gL_1g^{-1}=L_2$, $gL_2g^{-1}=L_1$ \cite[Lemma 2.5.8(i), Prop. 2.9.1(iv), Lemma 2.9.4]{KL}.

In case of $q=2$, $\Omega $ is not the only index $2$ subgroup of $S$, but we still have $\Omega ^+_4(2)=L_1L_2$, $L_i \lhd \Omega ^+_4(2)$, $i=1,2$, $L_1\cap L_2=1$, $L_i \cong SL_2(2)\cong S_3$, $i=1,2$.
Again $gL_1g^{-1}=L_2$, $gL_2g^{-1}=L_1$.
Let $N_i:=$ the unique normal subgrpup of $L_i$, $i=1,2$.
So $N_i$ is cyclic of order $3$, $i=1,2$, $|N_1N_2|=9$, and $gN_1g^{-1}=N_2$, $gN_2g^{-1}=N_1$.
\end{Note}

\begin{Lemma}\label{Thm2-L11}
Let $F=\bar{F}$, $\Char F=2$ and $S=SO^+_4(q)$, where $q$ is even.
Let $M$ be a non-trivial irreducible $FS$-module.
Then either:
\begin{enumerate}
\item $\dim_F M\ge 4$, or
\item $M$ is a trivial $\Omega$-module in case $q>2$, or a trivial $N_1N_2$-module in case $q=2$.
\end{enumerate}
\end{Lemma}

\textbf{Proof:}
Suppose by negation that $\dim_F M\le 3$.
Assume firstly that $M$ is a faithful $S$-module.
Recall that $S$ is generated by elements $x$ of order $2$.
Since $\Char F=2$, $(x-1)^2=0$ and $\ker(x-1) \supseteq \im(x-1)$.
Since $\dim \ker (x-1)+\dim \im(x-1)=\dim_FM \le 3$, it follows that each such $x$ is a transvection on $M$.

If $S$ acts primitively on $M$, the pair $(M,S)$ should appear in the list of \cite[Thm. 1.5]{KM}, but it does not.
If $S$ acts imprimitively on $M$ this implies, if $\dim_FM=3$, that $S \cong D \rtimes S_3$,  $D\cong \mathbb{Z}/m\mathbb{Z}\times \mathbb{Z}/m\mathbb{Z}$, $m$ is odd.
If $\dim_FM=2$, $S\cong D\rtimes S_2$, where $D\cong \mathbb{Z}/m\mathbb{Z}$, $m$ is odd. Both possibilities are inconsistent with the structure of $SO^+_4(q)$, $q$ even (if $q>2$, $SO_4^{+}(q)$ is not solvable, and $|SO_4^{+}(2)|=72)$.
Therefore $1 \neq N:=\ker(S\rightarrow SL(M))$ .

If $N\cap \Omega = 1$, then $N\times \Omega\approx S$, so $|N|=2$.
Hence $N$ is central in $S$ in contradiction to \cite[Tab. 2.1.D]{KL} (going to the natural $4$-dimensional $S$-module).
Therefore $1\neq K:=N\cap \Omega$.

If $K\cap L_1=1$, then $K$ commutes with $L_1$ (elementwise), so $K\subseteq L_2$. Therefore $N_2 \subseteq K$ (if $q=2$), or $L_2=K$ (if $q>2$).
Consequently $N_1=gN_2g^{-1}\subseteq gKg^{-1} \subseteq gNg^{-1}=N \Rightarrow N_1\subseteq N \cap \Omega =K$, if $q=2$.
So $N_1N_2\subset N$ in this case.
Or if $q>2$ then $L_1=gL_2g^{-1}\subseteq gKg^{-1} \subseteq gNg^{-1}=N$, so $L_1 \subseteq N \cap \Omega =K$, so $\Omega = L_1L_2 \subseteq N$ in this case.
A similar argument applies if $K\cap L_1 \neq 1$.$\Box$

\begin{Note}\label{Thm2-N11-1}
\underline{The case of $SO_4^{-}(q)$, $q$ is even.}

Using Kleidman-Liebeck \cite{KL}, we have the following properties:
\begin{enumerate}
\item $P \Omega _4^{-} (q)$ is simple \cite[Tab. 5.1.A]{KL};
\item Set $\Omega := \Omega _4^{-}(q)$. $Z(\Omega )=\Omega \cap F_q^*=1$ \cite[Prop. 2.9.3]{KL}, and consequently $P \Omega = \Omega$;
\item $[S:\Omega ]=2$, where $S=SO_4^{-}(q)$ \cite[Tab. 2.1.C]{KL};
\item $S$ is generated by reflections on its natural module \cite[Prop. 2.5.6, and (2.5.11)]{KL}.
\end{enumerate}
\end{Note}

\begin{Proposition}\label{Thm2-P11-2}
Let $M$ be a non-trivial irreducible $FS$-module, where $S=SO_4^{-}(q)$, $F=\bar{F}$, $\Char F=2$ and $q$ is even.
Then $\dim_FM\ge 4$.
\end{Proposition}

\textbf{Proof:}
We  shall firstly show that $M$ is a faithful $S$-module.
Indeed suppose by negation that $1\neq N:= \ker(S \rightarrow GL(M))$.
Clearly $N \lhd S$.
Since $\Omega $ is simple it follows that either $N\cap \Omega = 1$ or $N \cap \Omega = \Omega$.

If $N\cap \Omega = 1$, then Note \ref{Thm2-N11-1}(3) implies that $S= N \times \Omega$ and $|N|=2$.
So $N$ is central in $S$, $N=<g>$, $g^2=1$.
Consequently, by Schur's Lemma, $g \in \End _{S}(V)=F_q1$, where $V$ is the natural $4$-dimensional $S$-module.
But $g^2 =1$, $\Char F_q=2$ and $g \in F_q1$, imply that $g=1$, a contradiction.

If $N\cap \Omega = \Omega $, then $M$ is an irreducible non-trivial $S/\Omega \cong \mathbb{Z} _2$-module.
But this is impossible since irreducible $F\mathbb{Z} _2 $ -modules are trivial $1$-dimensional and $M$ is not.

By Note \ref{Thm2-N11-1}(4) we have that $S$ is generated by elements $x$ with $x^2=1$.
Suppose by negation that $\dim_F M \le 3$.
Hence $(x-1)^2=0$ implies that $\ker (x-1) \supseteq \im (x-1)$.
Since $\dim _F \ker (x-1) +\dim \im _F (x-1)=\dim _F M \le 3$, it follows that $x$ is a transvection on $M$.
So $S$ is generated by transvections on $M$.

If $S$ acts primitively on $M$, then the pair $(S, M)$ would have appeared in the list of \cite[Thm. 1.5]{KM}.
But it does not.

If $S$ acts imprimitively on $M$, then $S \cong D \rtimes S_3 $ or $S\cong D \rtimes S_2$, where $D$ is commutative, so $S$ is solvable, a contradiction.$\Box$

\begin{Note}\label{Thm2-N11-3}
\underline{The cases of $SO_6^{\pm}(q)$, $q$ is even.}

We shall treat both cases simultaneously.

Using Kleidman-Liebeck \cite{KL}, we have the following properties:
\begin{enumerate}
\item $P \Omega _6^{-}(q)$, $P \Omega _6^{+}(q)$ are both simple \cite[Tab. 5.1.A]{KL};
\item Set $\Omega := \Omega _6^{-}(q)$ if $S=SO_6^{-}(q)$, and $\Omega := \Omega _6^{+}(q)$ if $S=SO_6^{+}(q)$.
We have $Z(\Omega )=\Omega \cap F_q^*=1$ \cite[Prop. 2.9.3]{KL}, and consequently $P \Omega = \Omega$ in both cases;
\item $[S:\Omega ]=2$, where $S=SO_6^{\pm}(q)$, $q$ is even \cite[Tab. 2.1.C]{KL};
\item An irreducible non-trivial $F\Omega$-module, $\Char F =2$,
is either of dimension $4$ (the spin module) or of dimension $\geq 6$. \cite[Prop. 5.4.11]{KL}.
\end{enumerate}
\end{Note}

\begin{Lemma} \label {Thm2-L11-4}
$S=SO_6^{\pm}(q)$ has (respectively) two subgroups $H_1 \cong SO_2^{+}(q)$, $H_2 \cong SO_4^{\pm}(q)$ with $h_1h_2=h_2h_1$, $\forall h_1\in H_1$, $h_2 \in H_2$.
\end{Lemma}

\textbf{Proof:}
Let $S:=SO_6^{-}(q)$, and $\beta =\{ e_1, e_2, f_1, f_2, x, y \}$ the standard basis of $V$, the natural $6$-dimensional $O_6^{-}(q)$-module over $F_q$ \cite[Prop. 2.5.3 (ii)]{KL}.
Set $W= \Span_{F_{q}}\{e_2, f_2, x, y\}$.
We have $W^{\bot}=F_qe_1+F_qf_1$, $W$, $W^{\bot}$ are non-degenerate with respect to the restriction of $\kappa$, the quadratic form of $V$, and $V=W \bot W^{\bot}$.

As in the paragraph after \cite[Lemma 2.1.5]{KL} the isometry group $I(W,\kappa )\cong O_4^{-}(q)$ can be regarded as a subgroup of $O_6^{-}(q)$, by acting as an identity on $W^{\bot}$.
This shows that $SO_4^{-}(q)\cong S(W,\kappa):=H_2$ is a subgroup of $SO_6^{-}(q)$.
Similarly, since $I(W^{+}, \kappa ) \cong O_2^{-}(q)$, we have that $SO_2^{-}(q)\cong S(W^{\bot},\kappa):=H_1$ is also a subgroup of $SO_6^{-}(q)$ with elements acting as the identity on $W$.
Consequently $H_1$, $H_2$ commute element-wise.

The argument for $H_1 \cong  SO_2^{+}(q)$, $H_2\cong SO_4^{+}(q)$, $S=SO_6^{+}(q)$, is the same starting with the standard basis $\beta =\{e_1,f_1,e_2,f_2, e_3,f_3\}$  \cite[Prop. 2.5.3(i)]{KL}, of $V$, the natural $6$-dimensional $O_6^{+}(q)$-module.$\Box$

\begin{Proposition}\label{Thm2-P11-5}
Let $M$ be a non-trivial irreducible $FS$-module, where $S=SO_6^{\pm}(q)$, $F=\bar{F}$, $\Char F=2$ and $q$ is even.
Then $\dim_FM\ge 6$.
\end{Proposition}

\textbf{Proof:}
Suppose firstly that $\dim_F M=5$.

By Clifford's theorem \cite[Thm. (50.5)]{CR} we have that $M=M_1 \oplus \cdots \oplus M_r$, where $M_i$ is an homogenous component of an irreducible $F \Omega $-module, $i=1,...,r$.
Moreover $\dim _F M_1 = \dim _F M_i$, $\forall i $ and $S$ acts transitively and imprimitively on $\{M_1,...,M_r\}$.

Consequently $r\dim _FM_1=5$, so $r=1$ or $r=5$.
If $r=1$ then $M=V_1\oplus \cdots \oplus V_s$, where $V_1 \cong V_i$ as irreducible $F \Omega $-modules.
Consequently $\dim _F M=5 =s\dim _F V_1$.
By \ref{Thm2-N11-3}(4) we get that $\dim _F V_1 =\dim _F V_i=1$, $\forall i$.
Consequently as $\Omega$ is simple by \ref{Thm2-N11-3}(1) and \ref{Thm2-N11-3}(3), it follows that $V_i$ is a trivial $\Omega $-module $\forall i$ and therefore $M$ is a trivial $\Omega $-module.
This implies that $M$ is irreducible $S/\Omega \cong \mathbb{Z} _2$-module.
But $\Char F=2$ implies that the latter is a $1$-dimensional trivial $\mathbb{Z}_2$-module, a contradiction.
If $r=5$ then $\dim _F M_i = 1$, $\forall i$, and again $M$ is a trivial $\Omega $-module, leading to the same contradiction.

Suppose now that $\dim _F M \le 4$.
We shall next show that $M$ is a faithful $S$-module.

Suppose otherwise $1 \ne N:= \ker (S \rightarrow GL(M))$.
If $N \cap \Omega =1$ then \ref{Thm2-N11-3}(3) implies $S=N \times \Omega$, $|N|=2$ and consequently $N \subseteq Z(S)$.
This cannot hold since if $V$ is the natural irreducible $6$-dimensional $S$-module, then $Z(S) \subseteq \End _S (V)=F_q 1$.
So $N \subset F_q1$.
But $|N|=2$ and $\Char F_q=2$ imply that $|N|=1$, a contradiction.
If $N\cap \Omega \ne 1$ then by the simplicity of $\Omega $ we have that $N \cap \Omega = \Omega $ and therefore $M$ is an irreducible $S/\Omega \cong \mathbb{Z}_2$-module.
As before, this leads to a contradiction.

We shall show now that $M$ is a reducible $H_2$-module.

Suppose by negation that $M$ is an irreducible $H_2$-module.
Then by Lemma \ref{Thm2-L11-4} and the faithful action of $S$ on $M$ we have $H_1 \subseteq \End _{H_2}(M)=F1$, where the last equality is by Schur's Lemma.
Now $H_1 \cong SO_2^{\pm}\cong D_{2 (q \mp 1)}$ \cite[Sec. 2.9, p. 43]{KL}, so if $q \ne 2$ then $H_1 \subset F1$ violates the non-commutativity of the dihedral group $H_1$.
In case of $H_1 \cong SO_2^{+}(2)$, $H_1 \cong D_2$, the cyclic group of order 2,  but then $\Char F=2$ and $H_1 \subset F1$ imply that $H_1=1$, a contradiction.

Let $M=M_r \supset M_{r-1} \supset \cdots \supset M_1 \supset M_0 =0$ be a decomposition series of $FH_2$-modules.
Thus $M_i/M_{i-1}$ is an irreducible $FH_2$-module, with $\dim_F M_i/M_{i-1} \le 3$, $\forall i$.
By Lemma \ref{Thm2-L11} and Proposition \ref{Thm2-P11-2}, $M_i/M_{i-1}$ is a trivial $K$-module, where $K \subseteq H_2$ and $K$ is not a $2$-group ($K=H_2$ in case $H_2=SO_4^{-}(q)$ and $K=\Omega $ or $K=N_1N_2$ if $H_2=SO_4^{+}(q)$, according to  $q \ne 2$ or $q=2$ respectively).
Consequently $(g-1)^{2^r}(M) = 0 $, $\forall g \in K$.
By the faithful action of $S$ on $M$, this implies $g^{2^r}=1$, $\forall g \in K$, a contradiction.$\Box$

In order to deal with case (c), that is $n\ge 6$ and is even, $p=2$, $G=S_{n+1}$ or $S_{n+2}$, we need the following result of A. Wagner.

\begin{Theorem}\label{Thm2-T12} \cite{Wa}
$S_m$, with $m>6$ have a unique faithful modular $2$ representation of least degree, this degree being $m-1$, or $m-2$ according as $m$ is odd or even.
\end{Theorem}

\begin{Corollary}\label{Thm2-C13}
Let $G$ be as in case (c). Suppose $M$ is a non-trivial irreducible $FG$-module. Then either:
\begin{enumerate}
\item $\dim_FM \ge n$, or
\item $M$ is a trivial $A_{n+1}$-module, if $G=S_{n+1}$, or a trivial $A_{n+2}$-module, if $G=S_{n+2}$.
\end{enumerate}
\end{Corollary}

\textbf{Proof:}
Suppose $\dim_FM <n$. Then, by Thm. \ref{Thm2-T12}, $M$ is an unfaithful $G$-module. Let $N:=\ker(G \rightarrow GL(M))$. Then $N$ is normal in $G$ and the result follows.$\Box$

\begin{Lemma}\label{Thm2-L14}
Let $(n,p)=(2,3)$, $G=SL_2(5)$ be as in case (d).
Then $M_1=Fx_1+Fx_2$ is an irreducible $G$-module and $\deg x_1=\deg x_2$, where the notations are as in the proof of Thm. \ref{Thm2-T1}.
\end{Lemma}

\textbf{Proof:}
If $a_1=\deg x_1<\deg x_2=a_2$, then $M_1=Fx_1$, $B_2=F[x_1]_{a_2}$, $M_2=B_2+Fx_2/B_2$.
So $M_1$, $M_2$ are trivial $G$-modules (since $G$ is generated by elements of order $3$).
Therefore $(g-1)(M_i)=0$, $i=1,2$, implying that $(g-1)^2(x_j)=0$, $j=1,2$ $\forall g \in G$.
Consequently $(g^3-1)(x_j)=(g-1)^3(x_j)=0$, $j=1,2$, $\forall g \in G$.
This shows, by assumption (4), that $g^3=1$, $\forall g \in G$, so $G$ is a $3$-groups, an obvious absurd.
So $M_1=Fx_1+Fx_2$, and a similar argument ensures that $M_1$ is irreducible $G$-module.
This settles Thm \ref{Thm2-T1} in case (d).$\Box$

\begin{Lemma}\label{Thm2-L15}
Let $G=3.A_6$, $(n,p)=(3,2)$ be as in case (e).
Let $M$ be a non-trivial $FG$-module with $\dim_FM=2$, where $\Char F=2$.
Then $M$ is a reducible $G$-mocule.
\end{Lemma}

\textbf{Proof:}
Assume by negation that $M$ is irreducible.
Suppose firstly that $M$ is unfaithful.
Let $N:=\ker (G\rightarrow GL(M))$.
Let $<h>$ be the normal subgroup of $G$ satisfying $|<h>|=3$ and $G/<h>\cong A_6$.
If $N\cap<h>=1$, then $N\cong N<h>/<h>$ is a non-trivial normal subgroup of $G/<h>\cong A_6$.
Hence $N<h>=G$ and $N \cong A_6$.
But $G$ is generated by order $2$ elements, so the same holds for $G/N\cong <h>$, in contradiction to $|<h>|=3$.
Therefore $<h>\subseteq N$.

If $<h>=N$, then $M$ is a faithful $G/<h>$-module. But $G/<h>\cong A_6$. So by \cite[Prop. 5.3.7(ii)]{KL} $R_2(A_6)=3$, a contradiction.

If $<h>\subset N$ (properly), then $N/<h>$ is a non-trivial normal subgroup of $G/<h>\cong A_6$, implying that $N=G$, so $M$ is a trivial $G$-module, violating the assumptions.

All in all $M$ is a $2$-dimensional faithful irreducible $G$-module.
Recall that $G$ is generated by elements of order $2$.
Since $\dim_FM=2$ and $\Char F=2$, each such element is a transvection on $M$.
If $G$ acts primitively on $M$, then the pair $(M,G)$ would have appeared in the list of \cite[Thm. 1.5]{KM}, but it does not.

So we must conclude that $G$ acts imprimitively, faithfully and irreducibly on $M$.
Therefore we get $G \cong D\rtimes S_2$, where $D$ is abelian, so $G$ is solvable, a contradiction.$\Box$

\begin{Corollary}\label{Thm2-C16}
Let $G=3.A_6$, $(n,p)=(3,2)$ be as in case (e).
Then Theorem \ref{Thm2-T1} holds.
\end{Corollary}

\textbf{Proof:}
We assume by negation that $\dim_FM_i\le 2$, $\forall i$.
Let $\{X_{ji}\}$ be the decomposition series of $M_i$ for each $i$.
Then by Lemma \ref{Thm2-L15} $X_{ji}/X_{j-1,i}$ is a trivial $1$-dimensional $G$-module $\forall i$.
Consequently $(g-1)^{2^s}(M_i)=0$, $\forall i$, $\forall g\in G$, implying that $(g^{2^t}-1)(x_j)=(g-1)^{2^t}(x_j)=0$, $j=1,2,3$, $\forall g\in G$.
Therefore by assumption $(4)$ $g^{2^t}=1$, $\forall g \in G$, and $G$ is a $2$-group, an obvious contradiction.
Hence $\deg x_1=\deg x_2=\deg x_3$ and $M_1=Fx_1+Fx_2+Fx_3$ is an irreduccible $G$-module.$\Box$

\underline{The proof of Theorem 2 in case (ii)}

So $G\subset GL(U)$ is a irreducible imprimitive group, generated by transvections on $U$, $\dim_FU=n$, $F$ is algebraically closed.
So $G \subset SL(U)$.
By \cite[1.8, 1.9]{ZS} $G$ is a monomial subgroup and $\Char F = 2$.
Therefore $G= D\rtimes S_n$, where $D=\{\diag(\xi ^{a_1},...,\xi ^{a_n})|\sum_{i=0}^n a_i\equiv 0 (\bmod m)\}$, and $|<\xi >|=m$, $\xi \in F$, is a primitive $m$-th root of unity.
Therefore $\Char F=2$ implies that $m$ is an odd integer.

\begin{Note}\label{Thm2-N17}
\begin{enumerate}
\item $G \ne S_n$, since $U$ is an irreducible $G$-module.
\item $D \cong <\xi>\times\cdots\times<\xi>$, ($(n-1)$ times), so $|D|=m^{n-1}$.
\end{enumerate}
\end{Note}

\begin{Lemma}\label{Thm2-L18}
Let $M$ be an $FG$-module with $\dim_F M=n-i$.
Then there exists a subgroup $K \subseteq D$ having the following properties:
\begin{enumerate}
\item $|K|\ge m^i$;
\item $M$ is a trivial $K$-module.
\end{enumerate}
\end{Lemma}

\textbf{Proof:}
Let $\varphi :G\rightarrow GL(M)$ be the natural homomorphism.
$G$ is generated by order $2$ elements, implying that $\varphi (G)\subseteq SL(M)$.

Now $g^m=1$, $\forall g \in D$, and $(m,2)=1$, imply that $\varphi (D)$ is a commutative group of semi-simple transformations on $M$, and is therefore simultaneously diagonalizable ($F$ is algebraically closed).
Therefore after a base change of $M$ we have:
$$\varphi(D)\subseteq \{\diag (\lambda _1,...,\lambda_{n-i})|\lambda_j \in F \text{, } j=1,...,n-i\}.$$
So $\det (\varphi (D))=1$ implies that $\prod_{j=1}^{n-i}\lambda _j=1$.
Also $g^m=1$ implies that $\lambda_j^m=1$, so $\lambda_j=\xi^{a_j}$, where $\xi \in F$ is the $m$-th primitive root of unity, $j=1,...,n-i$.

Consequently:
$$\varphi(D)\subseteq \{\diag(\xi^{a_1},...,\xi^{a_{n-i}})|\sum_{j=1}^{n-i}a_j\equiv 0(\bmod m)\}\cong <\xi>\times\cdots\times <\xi>\text{, }(n-i-1)\text{-times}.$$
Hence $|\varphi(D)|\le m^{n-i-1}$.
Therefore $|K|\ge m^i$, where $K:=\ker(\varphi _{|D}:D\rightarrow GL(M))$.$\Box$

\begin{Corollary}\label{Thm2-C19}
Suppose one of the following holds:
\begin{enumerate}
\item $k=2$, that is $\dim_FM_1+\dim_FM_2=n$;
\item $\dim_FM_i=n-l$, for some $i$, $\dim M_j=1$, $\forall j\neq i$.
\end{enumerate}
Then $M_s$ is a trivial $K$-module $\forall s$, where $K\subseteq D$ is a subgroup with $K \neq 1$.
\end{Corollary}

\textbf{Proof:}
Set $K_i:=\ker (\varphi _i : D \rightarrow GL(M_i))$ $\forall i$.
If (1) holds then by Lemma \ref{Thm2-L18} $|K_i|\ge m^{n-\dim_FM_{i}}$, $i=1,2$.
Hence $|K_1||K_2|\ge m^{2n-n}=m^n > m^{n-1}=|D|$.
Therefore $K:=K_1\cap K_2\neq 1$, will do.

Suppose $(2)$ holds.
$G$ is generated by order $2$ elements, so $M_j$ is a trivial $G$-module $\forall j \neq i$.
Now by Lemma \ref{Thm2-L18} $|K_i|\ge m^l$, so $K:=K_i$ acts trivially on $M_s$ $\forall s$.$\Box$

\underline{The proof of Theorem 2, $G=D\rtimes S_n$, $n$ is even, $n > 6$}

Recall that $S_n$ is a subgroup of $G$.
As before, we consider the sequence of $G$-modules $M_i$, $i=1,...,k$.
If $\dim_FM_i\le n-3$ then by Theorem \ref{Thm2-T12}, $M_i$ is an unfaithful $S_n$-module.
Consequently $M_i$ is a trivial $A_n$-module $\forall i$.
This implies as before that $(g-1)^{2^s}(x_j)=0$, $\forall g \in A_n$, $j=1,...,n$, for some fixed s.
Therefore by the faithful action assumption of $G$ on $A=F[x_1,...,x_n]$, we have $g^{2^s}=1$ $\forall g \in A_n$.
So $A_n$ is a $2$-group, an obvious contradiction.

If $\dim_FM_i=n-2$, and $\dim_FM_j=2$, for some $j$.
Then by Corollary \ref{Thm2-C19}(1) $\{M_i,M_j\}$ are trivial $K$-modules for some $K\subseteq D$, $K\neq 1$.

If $\dim_FM_i=n-2$, $\dim_F M_j=\dim_FM_l=1$, where $j\neq l$.
Then by Corollary \ref{Thm2-C19}(2) $\{M_i,M_j,M_l\}$ are trivial $K$-modules, for some $K\subseteq D$, $K\neq 1$.

If $\dim_FM_i=n-1$, $\dim_F M_j=1$, then again by Corollary \ref{Thm2-C19}(1) $\{M_i,M_j\}$ are trivial $K$-modules for some $K \subseteq D$, $K\neq 1$.

So in all of these cases $(g-1)^{2^t}(x_j)=0$, $j=1,...,n$, $\forall g\in K$.
So $K$ is a $2$-group, in contradiction to $|K|$ being odd.
So $\deg x_1=\cdots =\deg x_n$, and $M_1=Fx_1+\cdots +Fx_n$ is irreducible $G$-module.$\Box$

\underline{The case of $G=D\rtimes S_4$ (that is $n=4$)}

By negation, we need to consider the following possibilities:
\begin{enumerate}
  \item $\dim_F M_i=3$, $\dim _F M_j=1$;
  \item $\dim_F M_i=2$, $\dim_FM_j=2$, $i\neq j$;
  \item $\dim_FM_i=2$, $\dim_F M_j=\dim_FM_l=1$, $j\neq l$;
  \item $\dim _F M_i=1$, $i=1,2,3,4$.
\end{enumerate}

In all cases we are led by Corollary \ref{Thm2-C19} to the existence of $K \subseteq D$, $K\neq 1$ and $M_l$ is a trivial $K$-module $\forall l$.
This implies as before that $K$ is a $2$-group, a contradiction.

\underline{The case of $G=D\rtimes S_2 $(that is n=2)}

We either have $\dim_FM_1=2$ and so $\deg x_1=\deg x_2$, or $\dim_FM_1=\dim_FM_2=1$ and we argue as before.$\Box$

\underline{The case of $G=D\rtimes S_6$ (that is n=6)}

By negation, we need to consider the following possibilities:
\begin{enumerate}
  \item $\dim_F M_i=5$, $\dim _F M_j=1$;
  \item $\dim_F M_i=4$, $\dim_FM_j=2$;
  \item $\dim_FM_i=4$, $\dim_F M_j=\dim_FM_r=1$, $j\neq r$;
  \item $\dim_F M_i=3$, $\dim_F M_j=3$, $i\neq j$;
  \item $\dim_FM_i=3$, $\dim_FM_j=2$, $\dim_FM_r=1$;
  \item $\dim_F M_i=3$, $\dim_F M_j=\dim_FM_r=\dim_FM_l=1$, for $3$ distinct indices $j,r,l$;
  \item $\dim_FM_i=\dim_FM_j=\dim_FM_r=2$, for 3 distinct indices $i,j,r$;
  \item $\dim_FM_i=\dim_FM_j=2$, $i\neq j$, $\dim_FM_r=\dim_FM_l=1$, $r\neq l$;
  \item $\dim_FM_i=2$, $\dim_FM_j=\dim_FM_t=\dim_FM_l=\dim_FM_r=1$, for $4$ distinct indices $j,t,l,r$;
  \item $\dim_FM_i=1$, $\forall i=1,...,6$.
\end{enumerate}

The only cases which are not directly handled by Corollary \ref{Thm2-C19} are (5), (7) and (8).

In case (5) we have by Lemma \ref{Thm2-L18}, that $|K_i|\ge m^3$, $|K_j|\ge m^4$, hence since $|D|=m^5$, $|K_i\cap K_j|=\frac{|K_i||K_j|}{|K_iK_j|}\ge \frac{|K_i||K_j|}{|D|}\ge\frac{m^7}{m^5}=m^2$.
Since $M_r$ is a trivial $G$-module, it follows that $\{M_i,M_j,M_r\}$ are trivial $K_i\cap K_j$-modules.
So the rest of the argument follows as before.

In case (7) $|K_i|\ge m^4$, $|K_j|\ge m^4$, $|K_r|\ge m^4$.
But $|K_iK_j||K_i\cap K_j|=|K_i||K_j|$, so $|K_i\cap K_j|=\frac{|K_i||K_j|}{|K_iK_j|}\ge \frac{|K_i||K_j|}{|D|}\ge\frac{m^8}{m^5}=m^3$.
Consequently $|(K_i\cap K_j)\cap K_r|=\frac{|K_i\cap K_j||K_r|}{|(K_i\cap K_j)K_r|}\ge \frac{m^3.m^4}{|D|}=m^2$.
Therefore $\{M_i,M_j,M_r\}$ are trivial $(K_i\cap K_j \cap K_r)$-modules.
So the rest follows as before.

In case (8), $|K_i|\ge m^4$, $|K_j|\ge m^4$, so as above $|K_i\cap K_j|\ge m^3$.
Now $\{M_r,M_l\}$ are trivial $G$-modules.
Hence $\{M_i,M_j,M_r,M_l\}$ are trivial $K_i\cap K_j$-modules, so we argue as before.$\Box$

\underline{The proof of Theorem 2, $G=D\rtimes S_n$, $n$ is odd, $n > 6$}

As before we consider the sequence of $G$-modules $M_i$, $i=1,...,k$.

If $\dim_FM_i\le n-2$, $\forall i$, then by Theorem \ref{Thm2-T12}, $M_i$ is an unfaithful $S_n$-module $\forall i$.
This implies that $M_i$ is a trivial $A_n$-module for each $i$.
Consequently $(g-1)^{2^s}(x_j)=0$, $j=1,...,n$, $\forall g \in A_n$, implying that $A_n$ is a $2$-group, an obvious contradiction.

If $\dim_FM_i=n-1$, $\dim_FM_j=1$, then by Corollary \ref{Thm2-C19} $\{M_i,M_j\}$ are trivial $K$-modules, for $1 \neq K \subset D$, a subgroup of $D$.
Consequently $(g-1)^{2^s}(x_j)=0$, $j=1,...,n$, $\forall g\in K$, implying that $K$ is a $2$-group, a contradiction.

All in all, $\deg x_1=\cdots = \deg x_n$, $M_1=Fx_1 +\cdots +Fx_n$ is an irreducible $G$-module.$\Box$

\underline{The case of $G=D\rtimes S_5$ (that is $n=5$)}

By negation, we need to consider the following possibilities:
\begin{enumerate}
\item $\dim_FM_i=4$, $\dim_FM_j=1$;
\item $\dim_FM_i=3$, $\dim_FM_j=2$;
\item $\dim_FM_i=3$, $\dim_FM_j=1=\dim_FM_l=1$, $j\neq l$;
\item $\dim_FM_i=\dim_FM_j=2$, $i\neq j$, $\dim_FM_l=1$;
\item $\dim_FM_i=2$, $\dim_FM_j=1$, $\forall j \neq i$;
\item $\dim_FM_i=1$, $i=1,...,5$.
\end{enumerate}

The only case not covered by Corollary \ref{Thm2-C19} is (4).
Here we get, by Lemma \ref{Thm2-L18}, $|K_i|\ge m^3$, $|K_j|\ge m^3$, hence $|K_i\cap K_j|=\frac{|K_i||K_j|}{|K_iK_j|}\ge\frac{m^3m^3}{|D|}=\frac{m^6}{m^4}=m^2$.
Since $G$ acts trivially on $M_l$, it follows that $\{M_i,M_j,M_l\}$ are trivial $K_i\cap K_j$-modules.
Therefore $(g-1)^4(x_s)=0$, $s=1,...,5$, $\forall g\in K_i\cap K_j$.
Hence $K_i\cap K_j$ is a $2$-group. A contradiction. The conclusion follows as before.$\Box$

\underline{The case of $G=D\rtimes S_3$, (that is $n$=3)}

By negation, we have to consider the following possibilities:
\begin{enumerate}
\item $\dim_FM_i=2$, $\dim_FM_j=1$;
\item $\dim_FM_i=1$, $i=1,,2,3$.
\end{enumerate}
Both cases are handled by Corollary \ref{Thm2-C19}, leading to a contradiction.$\Box$

\section{\bf Isolated quotient singularities in prime characteristic}\label{Sec-ThmB}

We assume throughout the present section that $G \subset SL(V)$ is a finite group, $\dim_FV$ is finite, $F$ is a field with $\Char F = p>0$.

We recall the following.

\begin{Theorem}\label{ThmB-T1} \cite[Theorem 6.1.11]{Wo}, \cite[Theorem 6.3.1]{Step}.
There exists a list of $\{(W,H)\}$, where $\dim_{\mathbb{\mathbb{C}}}W$ is finite, $H\subset GL(W)$ is a finite group, such that each isolated quotient singularities over $\mathbb{C}$ is isomorphic to one of $S(W)^H$.
\end{Theorem}

We shall refer to the above list as the Zassenhaus-Vincent-Wolf list.

The main result in this section is the following.

\begin{theoremb}
Suppose $S(V)^G$ is an isolated singularity and $F$ is algebraically closed.
Then $\widehat{S(V)^G}\cong \widehat{S(W)^H}$, where $(p,|H|)=1$, $\dim_FW=\dim_FV$ and $(W,H)$ is a $\bmod p$ reduction of a direct sum of members in the Zassenhaus-Vincent-Wolf list.
\end{theoremb}

\begin{Proposition} \label{ThmB-P4}
Suppose $G\subset SL(V)$, $F$ is perfect, $\Char F=p>0$ and $S(V)^G$ is an isolated singularity.
Then:
\begin{enumerate}
  \item $S(V)^{T(G)}$ is a polynomial ring;
  \item $T(G)=<P|P$ is a $p$-sylov subgroup of $G>$, and consequently $(|G/T(G)|,p)=1$.
\end{enumerate}
\end{Proposition}
Here $T(G)$ is the subgroup generated by all transvections on $V$ (or $V^*$).

\textbf{Proof:}
Let $0 \ne f \in V^*$ and $U=Ff$.
By Corollary \ref{ThmA-C5.8}  $S(V)^{G_U}$ is a polynomial ring, hence $G_U$ is generated by transvections implying $G_U\subseteq T(G)$.
Hence $G_U=T(G)_U$.
Therefore by Theorem A $S(V)^{T(G)}$ is a polynomial ring.

Moreover, let $P$ be a $p$-sylov subgroup of $G$, then by a classical result  $\exists f \in V^*$, $f \ne 0 $, with $P \subseteq G_U$, $U:=Ff$.
Since $G_U=T(G)_U$, it follows that $P \subseteq T(G)$.
Hence $<P|P$ is a $p$-sylov subgroup of $G > \subseteq T(G)$.

The converse inclusion is trivial since every transvection is of order $p$.
So $(2)$ is verified.$\Box$

\begin{Proposition} \label{ThmB-P4.5}
Let $G \subset SL(V)$ be a finite group and $F$ is perfect.
Suppose $S(V)^G$ is an an isolated singularity.
Then:
\begin{enumerate}
  \item $S(V)^G$ has rational singularity;
  \item $S(V)^G$ has a non-commutative crepant resolution if it is Gorenstein.
\end{enumerate}
\end{Proposition}

\textbf{Proof:}
By Proposition \ref{ThmB-P4} $S(V)^{T(G)}$ is a polynomial ring and $(|G/T(G)|,p)=1$.
Therefore $S(V)^G=(S(V)^T)^{G/T(G)}$ is a direct summand of $S(V)^{T(G)}$ (using $\frac{1}{r}\sum_{h \in G/T(G)}h$, as the Reinolds operator, where $r=|G/T(G)|$).
Consequently by \cite[Thm. 2.1]{HH} $S(V)^G$ is $F$-regular.
In particular it is $F$-rational.
Consequently by \cite[Cor. 1.10]{Ks} $S(V)^G$ has rational singularity.
This settles item (1).

We now show that the skew group ring $S(V)^{T(G)}*G/T(G)$ is a non-commutative crepant resolution (if $S(V)^G$ is Gorenstein).
We firstly recall that the extension $S(V)^G \subset S(V)^{(T(G)}$ is unramified in codimension 1 \cite[Lemma p.05]{Se}.
Consequently as in \cite[Cor. 2.24(2)]{Br1}, using $S(V)^G=(S(V)^{T(G)})^{G/T(G)}$, we have $S(V)^{T(G)}*G/T(G)\cong End _{S(V)^G}S(V)^{T(G)}$.
Since $(|G/T(G)|,p)=1$ it follows that $gl.\dim [S(V)^{T(G)}*G/T(G)]$ is finite \cite[Thm. 5.6]{McR}.
Now being a Cohen-Macaulay $S(V)^G$-module (since $S(V)^{T(G)}$ is such), it follows that $S(V)^{T(G)}*G/T(G)$ is a homologically homogenous PI ring.
The reflexivity of $S(V)^{T(G)}$ as a $S(V)^G$-module follows from the normality of $S(V)^{T(G)}$.
This settles item (2).$\Box$

The next property will be needed in the sequel.

\begin{Definition}\label{ThmB-D4.6}
We say that $V$ is void of fixed points if $V^g=0$, $\forall g \in G$, $g \neq 1$.
\end{Definition}

\begin{Remark} \label{ThmB-R4.7}
$V$ is void of fixed points if and only if $V^*$ is void of fixed points.
This is so since if $V^g\neq 0$, for some $g \in G$, then $(V^*)^{g^{-1}} \neq 0$.
To verify this, observe that if $A$ is the matrix of $g$ on $V$ with respect to a basis $\{v_1,...,v_n\}$ of $V$, then $A^t$ is the matrix of $g^{-1}$ on $V^*$ with respect to the dual basis $\{f_1,...,f_n\}$.
Then $v:=\sum_{i=1}^{n}\xi _i v_i \in V^g$ implies $\sum _{i=1}^{n} \xi _i f_i \in (V^*)^{g^{-1}}$.
\end{Remark}

The following result appears in \cite[Lemma 2.4]{Step} with the additional assumption $(|G|,p)=1$.

\begin{Lemma} \label{ThmB-L4.8}
Suppose $F$ is perfect and $G$ is free from pseudo-reflections.
Then the following are equivalent:
\begin{enumerate}
\item $S(V)^G$ is an isolated singularity;
\item $V$ is void of fixed points.
\end{enumerate}
\end{Lemma}

\textbf{Proof:}
By Corollary \ref{ThmA-C5.8}, (1) is equivalent to: $S(V)^{G_U}$ is a polynomial ring for each $U \subset V^*$, $\dim _FU=1$.
Suppose (1) holds.
If $V^g\neq 0 $, for some $g \in G$, $g \neq 1$, then by Remark \ref{ThmB-R4.7} $(V^*)^{g^{-1}}\neq 0$, hence $\exists U \subset V^*$, $\dim _F U = 1$ with $g^{-1}\in G_U$, so $G_U\neq 1$.
Now $S(V)^{G_U}$ is a polynomial ring, so $G_U$ is generated by pseudo-reflection, in contradiction to the assumption.
Conversely if (2) holds on $V$, then (2) holds on $V^*$.
Hence $G_U=1$, for each $U \subset V^*$, implying that $S(V)^{G_U}=S(V)$ is a polynomial ring, hence (1) holds.$\Box$

The following is well known.

\begin{Lemma} \label{ThmB-L4.9}
Let $(A,\mathfrak{m})$ be a local ring which is the localization of a finitely generated algebra over a field.
Then the following are equivalent:
\begin{enumerate}
\item $A$ satisfies $R_i$ (respectively $S_i$);
\item $\hat{A}$, the $\mathfrak{m}$-adic completion, satisfies $R_i$ (respectively $S_i$).
\end{enumerate}
\end{Lemma}
Consequently $A$ is an isolated singularity (respectively normal) if and only if $\hat{A}$ is such.

\textbf{Proof:}
The implication $(2) \Rightarrow (1)$ is valid for every Noetherian local ring as follows from \cite[Thm. 23.9(i)]{Ma}.

To prove the converse direction, recall that $A$ is a $G$-ring \cite[Cor,p.259]{Ma}.
Hence the map $A \rightarrow \hat{A}$ is regular \cite[p.256]{Ma}.
Consequently for every prime  $\mathfrak{p} \in \spec A$ the fiber ring $(\hat{A})_{\mathfrak{p}}/\mathfrak{p}(\hat{A})_{\mathfrak{p}}$ is geometrically regular.
Hence $(\hat{A})_{\mathfrak{p}}/\mathfrak{p}(\hat{A})_{\mathfrak{p}}$ satisfies $R_i$ and $S_i$.
Hence by \cite[Thm. 23.9(ii),(iii)]{Ma} $\hat{A}$ satisfies $R_i$ (respectively $S_i$).$\Box$

We shall need the following version of H. Cartan's theorem \cite[Lemma 2.3]{Step}:

\begin{Theorem}\label{ThmB-T2}
Let $F[[x_1,...,x_n]]$ be a power series ring and $H\subseteq \Aut F[[x_1,...,x_n]]$ is a finite group with $(|H|,p)=1$.
Suppose $\mathfrak{m}:=(x_1,...,x_n)$ is $H$-stable.
Then:
\begin{enumerate}
\item There are parameters $\{y_1,...,y_n\}$ such that $F[[x_1,...,x_n]]=F[[y_1,...,y_n]]$ and $H$ acts linearly on $W:=Fy_1+\cdots +Fy_n$;
\item $F[[x_1,...,x_n]]^H\cong \widehat{S(W)^H}$.
\end{enumerate}
\end{Theorem}

\textbf{Proof:}
For (1) we refer to \cite[Lemma 2.3]{Step}.
Now the inclusion $\widehat{S(W)^H}=\widehat{F[y_1,...,y_n]^H} \subseteq F[[y_1,...,y_n]]^H=F[[x_1,...,x_n]]^H$ is clear.
Note that by Galois theory $[Q(F[[x_1,...,x_n]]):Q(F[[x_1,...,x_n]])^H]=|H|$.
So we only need to affirm that $[Q(\widehat{S(W)}):Q(\widehat{S(W)^H})]\le |H|$, since this implies that $Q(\widehat{S(W)^H})=Q(F[[x_1,...x_n]])^H=Q(F[[x_1,...,x_n]]^H)$, and then the normality of $\widehat{S(W)^H}$ (which follows from the normality of $S(W)^H$)) implies that $\widehat{S(W)^H}=F[[x_1,...,x_n]]^H$.

Set $\mathfrak{n}=S(W)^H_{+}$, $\mathfrak{m}=S(W)_{+}$.
Then $\mathfrak{m}$ is the unique maximal above $\mathfrak{n}$, hence $\widehat{S(W)} =\lim _{\leftarrow i} S(W)/\mathfrak{n}^iS(W)$ implying that $\widehat{S(W)}=S(W) \otimes_{S(W)^H}\widehat{S(W)^H}$.
Hence  $[Q(\widehat{S(W)}):Q(\widehat{S(W)^H}]\le [Q(S(W)):Q(S(W)^H)]=|H|$.$\Box$

In order to translate the Zassenhaus-Vincent-Wolf list from objects over $\mathbb{C}$ to ones over an algebraically closed field $F$, with $\Char F = p >0$, we need the following result (well known to finite group theorists).
This is also sketched in \cite[Theorem 3.13]{Step}

\begin{Proposition} \label{ThmB-P6}
Let $G$ be a finite groups with $(|G|,p)=1$ and $F$ is algebraically closed.
Then:
\begin{enumerate}
  \item There exists an isomorphism:
    \begin{flalign*}
    \psi:\{\text{Irreducible } \mathbb{C}G -modules\} \rightarrow \{\text{ Irreducible } FG-modules\};
    \end{flalign*}
 \item $M^g=0$, $\forall 1 \neq g \in G$ if and only if $\psi (M)^g=0$, $\forall 1 \neq g \in G$.
\end{enumerate}
\end{Proposition}

\textbf{Proof:}
Observe firstly that by \cite[Thm. 83.5]{CR}, since $(|G|,p)=1$,
\begin{flalign*}
& \text{The number of non-isomorphic irreducible }  FG \text{-modules } \\
& = \text{ The number of }p\text{-regular conjugacy classes}\\
& = \text{ The number of conjugacy classes}\\
& = \text{ The number of non-isomorphic irreducible } \mathbb{C}G \text{-modules}.
\end{flalign*}
So the cardinality of both sets appearing  in $(1)$ is the same.

Let $K:= \mathbb{Q}(\xi )$, $\xi $ is a primitive $m$-th root of unity, where $m$ is a prime number, $p\ne m \ge \exp (G)$.
Then by Brauer's theorem \cite[Thm. (41.1)]{CR}, $K$ is a splitting field for $G$.
So $M=M_0 \otimes _{K} \mathbb{C}$, $M_0 \subset M$, and $M_0$ is an irreducible $KG$-module.
Let $R:=\mathbb{Z} [\xi ]$, so $R$ is the ring of integers of $K$ \cite[Thm. 2.6]{Wash}.
Let $P$ be a prime ideal in $R$, with $P \cap \mathbb{Z} = p \mathbb{Z} $.
Then there exists a $R_PG$-module $N$, $N \subset M_0$, with $N\otimes _{R_P} K=M_0$ (e.g. \cite[Thm. (75.2)]{CR}).

Now $N$ is an indecomposable $R_PG$-module.
Otherwise $M_0=N\otimes _{R_P}K$ would have been decomposable as $KG$-module.
Let $\bar{N}:=N/PN$.
Set $\bar{K}:=R_P/P_P=\mathbb{F}_p(\bar{\xi })\subset F$, $\bar{\xi }$ a primitive $m$-th root of the unity in $F$ (we identify $<\bar{\xi} >$ with the group of $m$-th roots of unity in $F$).
Then by \cite[Cor. (83.7]{CR} $\bar{K}$ is a splitting field for $G$.

$\bar{K}G$ is a semi-simple $\bar{K}$-algebra since  $(p,|G|)=1$, hence $\bar{N}$ is a projective $\bar{K}G$-module.
Consequently by \cite[Thm. (77.1)]{CR} $N$ is a projective $R_PG$-module, implying that it is a direct summand of a free $R_PG$-module.

We shall next show that $\bar{N}$ is an irreducible $\bar{K}G$-module.
Let $\hat{R}_P:=$ the $P_P$-adic completion of $R_P$.
Since $N$ is an indecomposable $R_PG$-module, it follows from the proof of \cite[Thm. (76.29)]{CR} that $\hat{R}_PN$ is an indecomposable $\hat{R}_PG$-module.
Also it is easily seen that $\hat{R}_PN$ is a projective $\hat{R}_PG$-module, implying that it is a direct summand of a free $\hat{R}_PG$-module.

Let $\bar{K}G=(\bar{K}G) \epsilon _1\oplus \cdots \oplus (\bar{K}G)\epsilon _m$ be a decomposition of $\bar{K}G$ into a direct sum of irreducible $\bar{K}G$ modules, where $\{\epsilon _1,...,\epsilon _m\}$ are orthogonal idempotents, with $\bar{1}=\epsilon _1+\cdots +\epsilon _m$.
By \cite[Thm. (77.11)]{CR} there are orthogonal idempotents $\{e_1,...,e_m\}$ in $\hat{R}_PG$, with $1 = e_1+\cdots +e_m$ and $\bar{e}_i=\epsilon_i$, $i=1,...,m$.
Therefore $\hat{R}_PG=(\hat{R}_PG)e_1\oplus \cdots \oplus (\hat{R}_PG)e_m$, and $\overline{(\hat{R}_PG)e_i}:=(\hat{R}_PG)e_i/P(\hat{R}_PG)e_i=(\bar{K}G)\epsilon _i$, $i=1,...,m$.
Hence $(\hat{R}_PG)e_i$ is an indecomposable projective $\hat{R}_PG$, $\forall i=1,...,m$.

Since $\hat{R}_PN$ is an indecomposable projective $\hat{R}_PG$-module, which is a direct summand of a free $\hat{R}_PG$-module, it follows from the Krull-Schmidt theorem for $\hat{R}_PG$-modules \cite[Thm. (76.26)]{CR} that $\hat{R}_PN\cong (\hat{R}_PG)e_i$, for some $i$.
Therefore:
$$\bar{N}=N/PN=\hat{R}_PN/P\hat{R}_PN\cong \overline{(\hat{R}_PG)e_i} =(\bar{K}G)\epsilon _i \text{, is an irreducible }\bar{K}G \text{-module}.$$

Moreover since $\bar{K}$ is a splitting field, it follows that $\bar{N}$ is an absolutely irreducible $\bar{K}G$-module.
So we define: $\psi (M):=\bar{N}\otimes _{\bar{K}}F$.
That this map is well defined (e.g. independent of the choice in $P$) follows from \cite[Cor. (82.2)]{CR}.

We finally need to verify that $\psi $ is $1-1$ (the ontoness will follow from the beginning paragraph).

Let $T,U$ be two irreducible representations of $KG$ with characters $\alpha, \beta$ (respectively).
Let $\bar{T}=\bar{U}$ denote the induced $\bar{K}G$-representations, and $\tau $ their common character..

Let $\chi $ be the Brauer character associated to $\tau $, e.g. \cite[p.589]{CR}.
Then $\overline{\chi (x)} =\tau (x)$, $\forall x \in G$ ($(p,|G|)=1$, so every $x\in G$ is $p$-regular).
On the other hand, since $\alpha, \beta$ are coming from $R_pG$-representations, it follows by \cite[p.589]{CR} that $\chi =\alpha$, $\chi = \beta $.
Therefore by \cite[(30.12), (30.14), (30.15)]{CR} $T$ and $U$ are equivalent.
This settles item $(1)$.

To prove $(2)$, suppose that $\psi (M)^g\neq 0$, then $\psi (M)=\bar{N} \otimes _{\bar{K}} F$ implies that $\bar{N}^g\neq 0$.
Let $0 \neq \bar{x} \in \bar{N}^g$, and $x\in N$ is some preimage of $\bar{x}$ in $N$.
Then $z:=\frac{1}{r} \sum_{i=0}^{r-1} g^i(x)\in N$, where $r=|g|$.
Moreover $g(z)=z$.
So $z\in N^g \subseteq M^g$.

We only need to affirm that $z\neq 0$.
This follows from $\bar{z}= \overline{\frac{1}{r}\sum_{i=0}^{r-1}g^i(x)} = \frac{1}{r}\sum_{i=0}^{r-1}g^i(\bar{x})=\bar{x}\neq 0$.

The reverse implication is easier.
If $M^g \neq 0$, then $M_0^g \neq 0$.
So pick $0 \neq x \in M_0^g$.
We can find $x_1 \in N^g$ with $x_1 \notin P_PN$ (say $x=\pi ^ i x_1$, $x_1 \notin P_PN$, $(\pi)=P_P)$.
Then $\bar{x}_1\in \bar{N}^g \subseteq \psi (M)^g$.
This settles $(2)$.$\Box$

\underline{The proof of Theorem B}

By Proposition \ref{ThmB-P4} $S(V)^{T(G)}$ is a polynomial ring and $(|G/T(G)|, p)=1$.
Set $H:=G/T(G)$.
By Cartan's theorem (Theorem \ref{ThmB-T2}) $\widehat{S(V)^{T(G)}}= F[[y_1,...,y_n]]$ and $G/T(G)$ acts linearly on $W:=Fy_1+...+Fy_n$.
Consequently $\widehat{S(V)^G} \cong (\widehat{S(V)^{T(G)}})^{G/T(G)}=F[[y_1,...,y_n]]^{G/T(G)}\cong \widehat{S(W)^H}$.
By Lemma \ref{ThmB-L4.9} the isolated quotient singularity property of $S(V)^G$ is translated to $\widehat{S(W)^H}$ and then to $S(W)^H$.

We next show that $H$ is free from pseudo-reflections on $W$.
Let $P(H)$ be the normal subgroup of $H$ generated by pseudo-reflections on $W$.
Let $K$ be a normal subgroup of $G$ satisfying $K/T(G)=P(H)$.
Then since $(|H|,p)=1$, $S(W)^{P(H)}$ is a polynomial ring.

We have $S(V)^K=(S(V)^{T(G)})^{K/T(G)}$.
So again by Cartan's theorem $\widehat{S(V)^K} = \widehat{S(V)}^K = \widehat{S(W)}^{P(H)} = \widehat{S(W)^{P(H)}}$ is a regular ring.
Consequently $S(V)^K$ is a polynomial ring.
Therefore $K$ is generated by pseudo-reflections on $V$ implying, since $G \subset SL(V)$, that $K \subseteq T(G)$ and $P(H)=1$.

Consequently by Lemma \ref{ThmB-L4.8} $W$ is void of fixed points, that is $W^g=0$, $\forall g \in H$, $g \neq 1$.

Now since $(|H|,p)=1$, $W$ is a completely irreducible $FH$-module, so $W=W_1\oplus \cdots \oplus W_r$, where $W_i$ is an irreducible $FH$-module and $W_i^g=0$, $\forall g \in H$, $g \neq 1$, $i=1,...,r$.

Consequently by Proposition \ref{ThmB-P6} there exists a $\mathbb{C}H$-module $M:=M_1\oplus \cdots \oplus M_r$, where $M_i$ is an irreducible $\mathbb{C}H$-module with $M_i^g=0$, $\forall g \in H$, $g \neq 1$, such that $\psi (M_i)=W_i$, $i=1,...,r$.
Now $M_i$, $i=1,...,r$ is a member of the Zassenhaus-Vincent-Wolf list (with respect to $H$) since $M_i^g=0$ $\forall g \in H$, $g \neq 1$  \cite[Thm. 5.3.1, 6.1.11, 6.3.1, 6.3.6]{Wo}.$\Box$

We shall now present modular isolated quotient singularity examples in arbitrary dimension ($\ge 3$). They are always Cohen-Macaulay, but not necessarily Gorenstein.

\begin{Example} \label{ThmB-E7}
Let $V=Fv+Fw_n \cdots +Fw_1$, be a $n+1$-dimensional $F$-vectorspace, with $n\ge 2$, and $a \in F$ is a primitive $l$-th root of unity.
We also require $(n,l)=1$.

Let $g:=\diag(a,...,a,a^{-n},a)$, $h:=1_{n+1}+e_{1,n+1}$.

The action on the basis elements is given by:

 $g(v)=av$,
  $\mbox{$g(w_i)$} = \left\{ \begin{array}{ll}
          \mbox{$aw_i$,} & \mbox{$i \ne 2$};\\
           \mbox{$ a^{-n}w_2$,} & \mbox{$i=2$}.
          \end{array} \right. $,
  $h(w_i)=w_i$, $i=1,...,n$, $h(v)=v+w_1$.

It is easy to verify that $h$ is a transvection, $g,h \in SL(V)$ and $gh=hg$.
Hence $G:=<g,h>$ is a commutative finite group, $|G|=pl$ and $G=\{g^ih^j|0\le i \le l-1$, $0\le j \le p-1\}$.
Set $H:=<h>$.

We next verify that $S(V)^G$ is an isolated singularity.
For that we show that $S(V)^{G_U}$ is a polynomial ring for each $U \subset V^*$ with $\dim_FU=1$.
Now by \cite[Lemma 2.1]{Br} $x \in G_U$ if and only if $(x-1)(V) \subseteq U^{\bot}$, where $U^{\bot}=\{y \in V|f(y)=0\}$, $U=Ff$.
Clearly $\dim_FU^{\bot}=\dim_FV-1=n$.

\underline{Case 1:} $w_r \notin U^{\bot}$, for some $r \in \{1,...,n\}$.

Let $x=g^ih^j\in G_U$, $0 \le i \le l-1$, $0\le j \le p-1$, then $(x-1)(V)\subseteq U^{\bot}$, hence:

$(x-1)(w_r)=(g^ih^j)(w_r)-w_r=g^i(w_r)-w_r$
  $\mbox{} = \left\{ \begin{array}{ll}
          \mbox{$(a^i-1)w_r$,} & \mbox{if $r \ne 2$};\\
           \mbox{$ (a^{-in}-1)w_2$,} & \mbox{if $r=2$}.
          \end{array} \right. $.

  Since $(l,n)=1$, it follows that $i=0$ and so $x = h^j \in H$ and $G_U \subseteq H$.
  Therefore either $G_U=1$ or $G_U=H$ and in both cases $S(V)^{G_U}$ is a polynomial ring.

\underline{Case 2:} $U^{\bot}=Fw_n+\cdots Fw_1$.

So $v \notin U^{\bot}$.
Similarly if $x=g^ih^j\in G_U$, then $(x-1)(v) \in U^{\bot}$.
But $(x-1)(v)=(g^ih^j)(v)-v=g^i(v+jw_1)-v=a^iv+ja^iw_1-v=(a^i-1)(v)+ja^iw_1$.
But $w_1\in U^{\bot}$, hence $(a^i-1)v\in U^{\bot}$, forcing $i=0$ and $x \in H$.
Hence $G_U \subseteq H$ and again $G_U=1$ or $G_U=H$ and in both cases $S(V)^{G_U}$ is a polynomial ring.

All isolated quotient singularities are Cohen-Macaulay as follows from \cite[Thm. 3.1]{Ke}.

Now $S(V)^H=F[v^p-w_1^{p-1}v, w_n,...,w_1]=S(M)$, where $M:=F(v^p-w_1^{p-1}v)+Fw_n+\cdots +Fw_1$.
We compute the action of $g^i$ on $M$:

$g^i(v^p-w_1^{p-1}v)= (g^i(v))^p-(g^i(w_1)^{p-1}g^i(v)=(a^iv)^p-(a^iw_1)^{p-1}a^iv= a^{ip}(v^p-w_1^{p-1}v)$, $g^i(w_r)$
 $\mbox{} = \left\{ \begin{array}{ll}
          \mbox{$a^iw_r$,} & \mbox{$r \ne 2$};\\
           \mbox{$ a^{-in}w_2$,} & \mbox{$r=2$}.
          \end{array} \right. $.

 Hence $S(V)^G=S(M)^{<g>}$ is Gorenstein if and only if $\det _Mg^i=1$, $\forall i$, if and only if $\deg_Mg=1$, that is $a^{p-1}=1$, so $a \in \mathbb{F}_p$.
 So if $a \in F -\{0\}$, then $S(V)^G$ is Gorenstein if and only if $l | p-1$.
\end{Example}


\section{\bf Acknowledgements}

Thanks are due to Yuval Ginosar, Martin Liebeck and Gunter Malle for their help in navigating through the turbulent waters of finite group theory.


\end{document}